\documentclass[a4paper,11pt]{amsart}
\usepackage{amssymb, graphicx}
\newtheorem{lemma}{Lemma}
\newtheorem{prop}{Proposition}
\newtheorem{thm}{Theorem}
\newtheorem{cor}{Corollary}
\theoremstyle{definition}
\newtheorem{rem}{Remark}
\newtheorem{defn}{Definition}

\newtheorem{claim}{Claim}
\newtheorem{assumption}{Assumption}
\newcounter{numl}

\newcommand{\labelnuml}{\textup{(\roman{numl})}}
\newenvironment{numlist}{\begin{list}{\labelnuml}%
{\usecounter{numl}\setlength{\leftmargin}{0pt}%
\setlength{\itemindent}{2\parindent}%
\setlength{\itemsep}{\smallskipamount}\def
\makelabel ##1{\hss \llap {\upshape ##1}}}}{\end{list}}
\newenvironment{bulletlist}{\begin{list}{\labelitemi}%
{\setlength{\leftmargin}{0pt}\setlength{\itemindent}{\parindent}%
\setlength{\itemsep}{\smallskipamount}\def
\makelabel ##1{\hss \llap {\upshape ##1}}}}{\end{list}}
\newcommand{\R}{{\mathbb R}}
\newcommand{\C}{{\mathbb C}}
\newcommand{\Z}{{\mathbb Z}}

\newcommand{\T}{{\mathbb T}}
\newcommand{\cA}{{\mathcal A}}
\newcommand{\cC}{{\mathcal C}}
\newcommand{\cL}{{\mathcal L}}
\newcommand{\cK}{{\mathcal K}}
\newcommand{\cM}{{\mathcal M}}
\newcommand{\cO}{{\mathcal O}}
\newcommand{\vE}{E}
\newcommand{\zE}{e}
\newcommand{\Lie}{{\mathcal L}}
\newcommand{\mult}{^{\scriptscriptstyle\times}}

\newcommand{\Scal}{\mathit{Scal}}
\newcommand{\Id}{\mathit{Id}}

\newcommand{\trace}{\mathop{\mathrm{tr}}\nolimits}
\newcommand{\pfaff}{\mathop{\mathrm{pf}}\nolimits}
\newcommand{\grad}{\mathop{\mathrm{grad}}\nolimits}
\newcommand{\vspan}{\mathop{\mathrm{span}}\nolimits}

\newcommand{\Aut}{\mathrm{Aut}}

\newcommand{\Mp}{p}
\newcommand{\Mpc}{p_{\mathrm{c}}}

\newcommand{\cb}{\overline{c}}
\begin{document}
\title[Hamiltonian 2-forms in K{\smash{\"a}}hler geometry, IV]
{Hamiltonian 2-forms in K{\smash{\"a}}hler geometry,\\
IV Weakly Bochner-flat K{\smash{\"a}}hler manifolds}
\author[V. Apostolov]{Vestislav Apostolov}
\address{Vestislav Apostolov \\ D{\'e}partement de Math{\'e}matiques\\
UQAM\\ C.P. 8888 \\ Succ. Centre-ville \\ Montr{\'e}al (Qu{\'e}bec) \\
H3C 3P8 \\ Canada}
\email{apostolo@math.uqam.ca}
\author[D. Calderbank]{David M.~J.~Calderbank}
\address{David M. J. Calderbank \\ Department of Mathematics\\
University of York\\ Heslington\\ York YO10 5DD\\ England}
\email{dc511@york.ac.uk}
\author[P. Gauduchon]{Paul Gauduchon}
\address{Paul Gauduchon \\ Centre de Math\'ematiques\\
{\'E}cole Polytechnique \\ UMR 7640 du CNRS
\\ 91128 Palaiseau \\ France}
\email{pg@math.polytechnique.fr}
\author[C. T\o nnesen-Friedman]{Christina W.~T\o nnesen-Friedman}
\address{Christina W. T\o nnesen-Friedman\\
Department of Mathematics\\ Union College\\
Schenectady\\  New York  12308\\ USA }
\email{tonnesec@union.edu}
\thanks{The first author was supported by NSERC grant OGP0023879, the second
author by an EPSRC Advanced Research Fellowship and the fourth author by the
Union College Faculty Research Fund.}
\date{\today}
\begin{abstract} We study the construction and classification of
weakly Bochner-flat (WBF) metrics (i.e., K\"ahler metrics with coclosed
Bochner tensor) on compact complex manifolds. A K\"ahler metric is WBF if and
only if its `normalized' Ricci form is a hamiltonian $2$-form: such $2$-forms
were introduced and studied in previous papers in the series. It follows
that WBF K\"ahler metrics are extremal.

We construct many new examples of WBF metrics on projective bundles and obtain
a classification of compact WBF K\"ahler $6$-manifolds, extending work by the
first three authors on weakly selfdual K\"ahler $4$-manifolds. The
constructions are independent of previous papers in the series, but the
classification relies on the classification of compact K\"ahler manifolds with
a hamiltonian $2$-form~\cite{ACGT3}.
\end{abstract}
\maketitle

A K\"ahler metric is said to be \emph{weakly Bochner-flat} (WBF) if the
Bochner tensor (a component of the curvature tensor) is coclosed. By the
differential Bianchi identity, this is equivalent to an overdetermined first
order linear equation on the Ricci form $\rho$. Examples include Bochner-flat
K\"ahler metrics (where the Bochner tensor is zero)---in particular metrics of
constant holomorphic sectional curvature (CHSC)---and products of
K\"ahler--Einstein metrics (for which $\rho$ is parallel).

The equation satisfied by the Ricci form of a WBF K\"ahler metric means that
the \emph{normalized Ricci form} $\tilde\rho:= \rho - \frac{\Scal_g}{2 (m +
1)} \omega$ is a \emph{hamiltonian $2$-form}: a real $(1,1)$-form (i.e., a
$J$-invariant $2$-form) $\phi$ on a K\"ahler manifold $(M,J,g,\omega)$, of
real dimension $2m>2$ is said to be \emph{hamiltonian}~\cite{ACG2} if
\begin{equation*}
2 \nabla _{X} \phi = d \trace\phi \wedge (JX)^{\flat} - (J d \trace\phi ) 
\wedge X^{\flat}
\end{equation*}
for all $X \in TM$ (where $X^{\flat}(Y) = g(X,Y)$ for $Y \in TM$ and
$\trace\phi = \langle\omega, \phi\rangle_{g}$).

The \emph{momentum polynomial} of a hamiltonian $2$-form $\phi$ is
\begin{equation*}
\Mp(t) := (-1)^m\pfaff (\phi-t\omega)
=  t^m - (\trace\phi) \,t^{m-1} + \cdots + (-1)^m \pfaff\phi,
\end{equation*}
where the \emph{pfaffian} is defined by $\phi \wedge \cdots \wedge \phi =
(\pfaff\phi)\omega\wedge\cdots\wedge\omega$. The reason for calling $\phi$
\emph{hamiltonian} is that the functions $\Mp(t)$ on $M$ (for $t\in\R$) are
Poisson-commuting hamiltonians for Killing vector fields $K(t) := J \grad_g
\Mp(t)$~\cite{ACG2}. The integer $\ell=\max_{x\in M} \dim\vspan
\{K(t)_x:t\in\R\}$ is called the \emph{order} of the hamiltonian $2$-form (and
$0\leq \ell\leq m$). The order of a WBF metric is defined to be the
order of its normalized Ricci form. Note that the Fubini--Study metric on $\C
P^m$ has order zero, but admits hamiltonian $2$-forms of any order $0\leq
\ell\leq m$~\cite{ACG2}.

It follows that WBF K\"ahler metrics are extremal in the sense
of~\cite{calabi1}. We thus have the following implications between classes of
K\"ahler metrics:
\begin{equation}\label{inclusions}\begin{tabular}{ccccc}
CHSC & $\Rightarrow$ &K\"ahler--Einstein & $\Rightarrow$ & CSC\\
$\Downarrow$ & &$\Downarrow$ & &$\Downarrow$\\
Bochner-flat &$\Rightarrow$& WBF &$\Rightarrow$&Extremal.
\end{tabular}
\end{equation}

The observation that a K\"ahler metric is WBF if and only if the normalized
Ricci form is hamiltonian motivated us to indulge in a detailed study of the
local and global theory of hamiltonian $2$-forms on K\"ahler
manifolds~\cite{ACG2,ACGT3} as well as the application of this to the theory
of extremal K\"ahler metrics~\cite{ACGT4}. For the final paper in this series,
we are now returning to our initial interest in WBF K\"ahler metrics.

We do not wish to impose the study of hamiltonian $2$-forms on the reader of
this paper, so we therefore propose to make the \emph{constructions} of WBF
metrics herein essentially self-contained, whereas for the \emph{necessity} of
the form of these constructions (both as motivation and as the source of the
classification results we obtain) we review in section~\ref{preliminaries} the
facts we require from the general theory. These results will allow us to
classify WBF metrics on compact $6$-manifolds.

The structure of the paper is as follows. In section~\ref{preliminaries} we
review the general theory of K\"ahler metrics with hamiltonian
$2$-forms~\cite{ACG2,ACGT3,ACGT4} with a special attention to the case when
the hamiltonian form has order $\ell=1$. We present an explicit construction
of such metrics on a class of `admissible' projective bundles of the form
$M=P(\vE_0\oplus\vE_\infty)\to S$, where $\vE_0$ and $\vE_\infty$ are
projectively flat hermitian vector bundles over a K\"ahler manifold $S$
endowed with compatible local product structure. According
to~\cite{ACGT3,ACGT4}, any K\"ahler manifold admitting a hamiltonian $2$-form
of order $1$ is obtained by this construction up to a covering, and if there
is no torsion in $H^2(S,\cO)$, we can take the covering to be trivial.

In section~\ref{KEsec}, as a warm-up, we use K\"ahler--Ricci
solitons~\cite{koiso} to study K\"ahler--Einstein metrics on admissible
bundles $M=P(\vE_0\oplus\vE_\infty)\to S$ with $S$ is a product of positive
K\"ahler--Einstein manifolds. We show that a K\"ahler--Ricci soliton exists
(and is unique) if and only if $M$ is a Fano manifold.  These examples were
found by Koiso~\cite{koiso}, and the vanishing of the Futaki invariant is
necessary and sufficient for the existence of a K\"ahler--Einstein metric,
cf.~\cite{koiso}.

In the remainder of the paper, we study WBF metrics in general.  In
section~\ref{WBFsec} we construct many compact WBF manifolds of order $1$,
including all such examples in dimension $6$.  This leads to a classification
of WBF $6$-manifolds $M$ in section~\ref{6-mnf-WBF}: they are either order $0$
and generalized K\"ahler--Einstein, or they are order $1$, and---apart from
one example on $P(\cO\oplus\cO(1)\otimes\C^2)\to \C P^1$---are then projective
line bundles over a ruled surface or a positive K\"ahler--Einstein surface. In
each case the WBF K\"ahler metric is unique up to scale and pullback by an
automorphism of $(M,J)$.

This is much richer than the classification of WBF $4$-manifolds, where the
only example of order $1$ is the first Hirzebruch surface
$P(\cO\oplus\cO(1))\to \C P^1$~\cite{ACG1}. It is natural to conjecture that
all compact WBF K\"ahler manifolds have order $0$ or $1$, but such a result is
out of reach using the explicit methods of this paper.

\tableofcontents

\section{Hamiltonian 2-forms and WBF K\"ahler metrics}\label{preliminaries}

We begin by recalling the classification of compact K\"ahler manifolds with a
hamiltonian $2$-form from~\cite{ACG2, ACGT3, ACGT4}, focussing on the case
that the hamiltonian $2$-form has order $1$. The output of this classification
is a self-contained Ansatz that we shall use to construct WBF K\"ahler metrics
in section~\ref{WBFsec}, so that we only need the results of~\cite{ACG2,
ACGT3, ACGT4} for the classification results we obtain. We adopt the notations
and conventions of~\cite{ACGT4} and refer to~\cite[\S1 \& App. A]{ACGT4} for
further information.

\subsection{Classification of hamiltonian $2$-forms}

Let $(M,g,J,\omega)$ be a compact connected K\"ahler $2m$-manifold with a
hamiltonian $2$-form $\phi$ of order $\ell$. Then, according to~\cite{ACGT3}
the vector fields $\{K(t):t\in\R\}$ described in the introduction generate an
effective isometric hamiltonian action of an $\ell$-torus $\T$ on $M$. The
stable quotient $\smash{\hat S}$ of $M$ by the induced action of the
complexified torus $\T^c$ is covered by a product of K\"ahler manifolds $S_a$
indexed by the distinct constant roots of $\Mp(t)$, the dimension of $S_a$
being $2d_a$, where $d_a$ is the multiplicity of the corresponding root.

It was also shown in~\cite{ACGT3,ACGT4} that there is a subset $\cA$ of the
constant roots such that $M$ is a projective bundle over a complex manifold
$S$ covered by $\prod_{a\in\cA} S_a$ in such a way that $\smash{\hat S}$ is a
fibre product of flat projective unitary bundles over $S$, indexed by the
remaining constant roots. In this paper, we shall always be in a situation
where the following assumption holds for these bundles.

\begin{assumption}\label{assume-P(E)} A flat projective unitary
$\C P^r$-bundle on $S$ is of the form $P(E)$, where $E$ is a rank $r+1$
projectively-flat hermitian holomorphic vector bundle.
\end{assumption}
If $S$ is simply connected, then any flat projective unitary $\C P^r$-bundle
is trivial, hence of the form $P(E)$ with $E\cong{\mathcal E}\otimes \C^{r+1}$
for a holomorphic line bundle ${\mathcal E}$.  In general the obstruction to
the existence of $E$ is given by a torsion element of $H^2(S, \cO^*)$
(cf.~\cite{elencwajg-narasimhan}). In particular, such an $E$ always exists if
$S$ is a Riemann surface.

It then follows, as in~\cite[App. A]{ACGT4}, that by formally adjoining
additional constant roots of multiplicity $0$ (corresponding to $\C P^0$
bundles over $S$) that we can write $\smash{\hat S}=P(\vE_0)\times_S
P(\vE_1)\times_S\cdots\times_S P(\vE_\ell)\to S$, where $\vE_j\to S$ are
projectively-flat hermitian bundles of ranks $d_j+1$ $(d_j\geq 0)$, which can
be chosen so that $M=P(\vE_0\oplus\vE_1\cdots\oplus\vE_\ell)\to S$. Thus the
distinct constant roots are labelled by
$\smash{\hat\cA}:=\cA\cup\{0,1,\ldots\ell\}$, and $S_a\cong\C P^{d_a}$ for
$a\in\{0,1,\ldots\ell\}$.  We remark that $M$ has a blow-up of the form
$\smash{\hat M}=P(\cL_0\oplus\cL_1\cdots\oplus\cL_\ell)\to \smash{\hat S}$ for
line bundles $\cL_j$.  If $d_j=0$ for all $j\in\{0,1,\ldots \ell\}$ then
$\smash{\hat M}=M$ and $\smash{\hat S}=S$. Otherwise we say {\it a blow-down
occurs}.

The extreme cases $\ell=0$ and $\ell=m$ are quite straightforward.
\begin{bulletlist}
\item If $\ell=0$, $M=\smash{\hat S}=S$ is a local K\"ahler product and the
hamiltonian $2$-form $\phi$ is a constant linear combination of the
corresponding K\"ahler forms.
\item If $\ell=m$, $(M,J)$ is biholomorphic to $\C P^m$ (and $\smash{\hat
  S}=S$ is a point).
\end{bulletlist}

For the intermediate cases, there is also an explicit description, but we
shall only need it in the case $\ell=1$ to which we now turn. Here it is
convenient to index the constant roots by
$\smash{\hat\cA}=\cA\cup\{0,\infty\}$ so that $\cA$ can be taken as a finite
subset of $\Z^+$.

\subsection{Admissible bundles and metrics}\label{projbundle}

\begin{defn}
A projective bundle of the form $M = P\bigl(\vE_0 \oplus \vE_\infty \bigr)
\stackrel{p}{\to} S$ will be called \emph{admissible} or an \emph{admissible
manifold} if:
\begin{bulletlist}
\item $S$ is a covered by a product $\smash{\tilde S}=\prod_{a\in \cA}
S_a$ (for $\cA\subset\Z^+$) of simply-connected K\"ahler manifolds
$(S_a,\pm g_a,\pm \omega_a)$ of real dimensions $2d_a$;
\item $\vE_0$ and $\vE_{\infty}$ are holomorphic projectively-flat hermitian
vector bundles over $S$ of ranks $d_0+1$ and $d_\infty+1$ with
$\cb_1(\vE_\infty)-\cb_1(\vE_0) = [\omega_S/2\pi]$ and $\omega_S=\sum_{a\in
\cA}\omega_a$.
\end{bulletlist}
\end{defn}
\noindent In the first condition, it is convenient to let $(g_a,\omega_a)$ be
positive or negative definite: otherwise we would have to admit signs in the
definition of $\omega_S$.  The second condition means that we can fix
hermitian metrics on $E_0$ and $E_\infty$ whose Chern connections have
tracelike curvatures $\Omega_0\otimes\Id_{E_0}$ and
$\Omega_\infty\otimes\Id_{E_\infty}$ satisfying
$\Omega_\infty-\Omega_0=\omega_S$. We normalize the induced
fibrewise Fubini--Study metrics $(g_0,\omega_0)$ and
$(-g_\infty,-\omega_\infty)$ on $P(E_0)$ and $P(E_\infty)$ to have scalar
curvatures $2d_0(d_0+1)$ and $2d_\infty(d_\infty+1)$.

We also have $\smash{\hat M}=P(\cO\oplus\smash{\hat\cL})\to\smash{\hat S}$
with $c_1(\smash{\hat\cL}) = [\omega_{\hat S}/2\pi]$ and $\omega_{\hat
S}=\sum_{a\in{\hat\cA}} \omega_a$.

\begin{rem}\label{integrality}
The existence of the line bundle $\smash{\hat\cL}\to\smash{\hat S}$ with
$c_1(\smash{\hat \cL})=[\omega_{\hat S}/2\pi]$ implies that $\omega_{\hat S}$
is \emph{integral} in the sense that $[\omega_{\hat S}/2\pi]$ is in the image
of $H^2(\smash{\hat S},\Z)$ in $H^2(\smash{\hat S},\R)$. When $S$ is a global
K\"ahler product (so we have $M=P(\cO\otimes\C^{d_0+1}\oplus\cL\otimes
\C^{d_\infty+1})\to S=\prod_{a\in\cA}S_a$) this integrality condition means
that each $\omega_a$ is integral, i.e., the compact manifolds $(S_a,\pm
g_a,\pm\omega_a)$ are Hodge. We write $\omega_a=q_a \alpha_a$ for an integer
$q_a\neq 0$, where $\alpha_a$ is a primitive integral K\"ahler form on $S_a$,
so that $q_a$ is a nonzero integer with the same sign as $(g_a,\omega_a)$, and
$q_0=1$ and $q_\infty=-1$.

If $\pm g_a$ is K\"ahler--Einstein, then $\rho_a=p_a\alpha_a$ and where $p_a$
is an integer (called the \emph{Fano index} for positive K\"ahler--Einstein
metrics).  We set $s_a=p_a/q_a$ and then $\Scal_a=\pm 2 d_a s_a$, where the
sign is that of $q_a$, so the scalar curvature of $\pm g_a$ has the same sign
as $p_a$.  For instance, if $S_{a}$ is $\C P^{1}$ and $g_{a}$ is negative
definite (i.e., $q_a$ is negative), then $\Scal_{a}$ is positive (and $p_a$ is
positive), but $s_a$ is negative. By the well-known Kobayashi--Ochiai
inequality~\cite{kob-och} $p_a\leq d_a+1$, where equality holds iff $S_a=\C
P^{d_a}$. Comparing the Chern classes $c_1(\cL_a)=[q_a \alpha_a/2\pi]$ and
$c_1(\cK^{-1})=[p_a\alpha_a/2\pi]$, we have that $\cL_a^{p_a}$ is $\cK^{-q_a}$
tensored by a flat line bundle.  If $p_a$ is not zero (i.e., $S_a$ is not
Ricci-flat), this gives $\cL_a\cong\cK^{-q_a/p_a}\otimes\cL_{a,0}$ for some
flat line bundle $\cL_{a,0}$. For instance if $S_a=\C P^{d_a}$, then
$p_a=d_a+1$ and $\cL_a\cong \cO(q_a)$.
\end{rem}

We now describe the K\"ahler metrics which admit a hamiltonian $2$-form $\phi$
of order $\ell=1$. In this case the hamiltonian torus action is just an $S^1$
action generated by a single hamiltonian Killing vector field $K=J\grad_g z$,
and without loss, we can take the image of its momentum map $z$ to be
$[-1,1]$. We denote the constant roots by $-1/x_a$ and we have that $0< |x_a|
\leq 1$ with equality iff $a\in\{0,\infty\}$; we can take $x_0=1$ and
$x_\infty=-1$.  Then $M^0:=z^{-1}((-1,1))$ is a principal $\C\mult$-bundle
over $\smash{\hat S}$ with connection $1$-form $\theta$
\textup($\theta(K)=1$\textup) and there are K\"ahler metrics $(\pm
g_a,\pm\omega_a)$, which are Fubini--Study metrics for $a\in\{0,\infty\}$,
with the signs chosen so that $\omega_a/x_a$ is positive for all $a$, together
with a smooth function $\Theta$ on $[-1,1]$ such that the K\"ahler structure
on $M^0$ is
\begin{equation}\label{metric}\begin{split}
g&=\sum_{a\in\smash{\hat\cA}\vphantom{I}} \frac{1 + x_a z}{x_a} g_a
+\frac{dz^{\smash 2}}{\Theta(z)}+\Theta(z)\theta^2,\\
\omega&=\sum_{a\in\smash{\hat\cA}\vphantom{I}} \frac{1 + x_a z}{x_a}\omega_a
+dz\wedge \theta,\qquad\text{where}\qquad
d\theta=\sum_{a\in\smash{\hat\cA}\vphantom{I}} \omega_a,
\end{split}\end{equation}
and $\Theta$ satisfies
\begin{gather}\label{positivity}
\Theta>0\quad\text{on}\quad (-1,1),\\
\Theta(\pm 1) = 0,\qquad
\Theta'(\pm 1) = \mp 2.
\label{boundary}
\end{gather}
It follows from~\cite{ACGT3,ACGT4} that if $M$ admits a hamiltonian $2$-form
of order $1$ and either Assumption~\ref{assume-P(E)} holds or no blow-downs
occur, then $M=P(\vE_0\oplus\vE_\infty)\to S$ is an admissible bundle, and the
above conditions are necessary and sufficient for the compactification of a
metric of the form~\eqref{metric} on $M$, where $z\colon M\to [-1,1]$ with
$P(E_0\oplus 0)=z^{-1}(1)$ and $P(0\oplus E_\infty)=z^{-1}(-1)$, $\theta$ is a
connection $1$-form (see~\cite{ACGT4} for more details), the $S^1$ action
generated by $K$ is given by scalar multiplication in $E_\infty$ (or
equivalently in $E_0$), and the local product structure in~\eqref{metric}
coincides with the given local product structure on $\smash{\hat
S}=P(E_0)\times_S P(E_\infty)\to S$.

We refer to a compatible metric of the form~\eqref{metric} on an admissible
bundle as an {\it admissible metric}. It is straightforward (and standard) to
see that the conditions~\eqref{positivity}--\eqref{boundary} are sufficient
for the compactification of metrics of the form, so that we can regard the
above as an Ansatz for constructing K\"ahler metrics on admissible bundles,
independently of the theory of hamiltonian $2$-forms.

\subsection{WBF K\"ahler metrics of order $0$ and $1$}

According to the theory of hamiltonian $2$-forms, a WBF K\"ahler manifold $M$
of order $0$ is a local K\"ahler product and the normalized Ricci form is a
constant linear combination of the corresponding K\"ahler forms. It follows
that $M$ is generalized K\"ahler--Einstein (i.e., its universal cover is a
product of K\"ahler--Einstein manifolds).

In the order $1$ case, we have the following characterization of WBF K\"ahler
metrics of the form~\eqref{metric}.
\begin{prop}\label{equations}
Let $(g,J,\omega)$ be a K\"ahler metric with a hamiltonian $2$-form $\phi$ of
order $1$ as in~\eqref{metric}, and write $F(t)=\Theta(t)p_c(t)$ with $p_c(t)
=\prod_{a\in\hat\cA} (1+ x_a t)^{d_a}$.  Then $g$ is WBF, with $\tilde\rho$ a
constant linear combination of $\phi$ and $\omega$, iff
\begin{itemize}
\item $F'(t)= Q(t) p_c(t)$ and $Q$ is a
polynomial of degree $\leq 2$\textup;
\item for all $a$, $\pm g_a$ is K\"ahler--Einstein with scalar curvature $\pm
d_a Q(-1/x_a)$.
\end{itemize}
$g$ is then K\"ahler--Einstein iff $Q$ has degree $\leq 1$.
\end{prop}
(Here we use the conventions of~\cite{ACGT4}, so that, compared to
\cite{ACG2}, we have $\eta_a=-1/x_a$ and have rescaled $F(z)$ and $p_c(z)$ by
$\prod_{a\in\hat\cA} x_a$.)

For the necessity of these conditions when $(g,J,\omega)$ is WBF, we refer
to~\cite{ACG2}, but their sufficiency is a straightforward verification.
Together with the discussion of the previous paragraph, we therefore have an
Ansatz for constructing admissible WBF K\"ahler metrics on admissible
projective bundles.

\section{K\"ahler--Einstein metrics and K\"ahler--Ricci solitons} \label{KEsec}

Recall that a {\it K\"ahler--Ricci soliton} on a compact complex manifold
$(M,J)$ is a compatible K\"ahler metric $(g,\omega)$ satisfying
\begin{equation}\label{rs1}
\rho - \lambda \omega = \Lie_V\omega,
\end{equation}
where $V$ is a real holomorphic vector field with zeros and $\lambda$ is a
real constant (necessarily equal to $\int_M \Scal_g \,\omega^m /\int_M
\omega^m$). It follows from \eqref{rs1} that the Futaki invariant
$\mathfrak{F}_{[\omega]}(V)$ vanishes iff the metric is K\"ahler--Einstein: if
$V = J \grad_g f + \grad_g h$, $\cL_V\omega=dd^c h$ and the imaginary part of
$\mathfrak{F}_{[\omega]}(V)$ reduces, after integrating by parts, to a nonzero
multiple of the $L^2$-norm of $\grad_g h$; if this is zero, $V$ is a
hamiltonian Killing vector field, so $\cL_V\omega=0$. Note that if $V$ is
nonzero then by the Bochner formula $\lambda>0$, and so $c_1(M)$ is positive,
i.e., $(M,J)$ is a Fano manifold.

The theory of K\"ahler--Ricci solitons on Fano manifolds has recently received
attention as a natural generalization of K\"ahler--Einstein metrics. In
particular, a number of uniqueness results for such metrics have been
established~\cite{tian-zhu1,tian-zhu2}, as well as existence results in the
case of toric Fano manifolds~\cite{wang-zhu} and certain geometrically ruled
complex manifolds~\cite{koiso}.

We now adapt arguments from \cite{koiso} to construct (admissible)
K\"ahler--Ricci solitons on admissible projective bundles $M= P(\cO 
\otimes \C^{d_{0}+1} \oplus \cL \otimes \C^{d_{\infty}+1}) \to S$, by 
taking $V=(c/2) \grad_g z$ for a real constant
$c$. Since $\Lie_V\omega = (c/2) dd^c z$ and
\begin{equation}\label{rho}
\rho =
\sum_{a} \rho_{a} - \frac{1}{2} d d^c \log F = \sum_{a} \rho_a - 
\frac{1}{2}\frac{F'(z)}{\Mpc(z)} \sum_{a} \omega_{a}
-\frac{1}{2}\Bigl(\frac{F'}{\Mpc}\Bigr)'(z) dz \wedge \theta,
\end{equation}
where $F$ and $p_c$ are as defined in Proposition \ref{equations}
(see~\cite{ACG2}), \eqref{rs1} is equivalent to
\begin{gather} \label{rs2}
\sum_{a} \rho_{a} = \sum_{a}\frac{1}{2}\biggl(\frac{F'(z)}{\Mpc(z)} + 
c\frac{F(z)}{\Mpc(z)} + 2\lambda(z + 1/x_{a})\biggr) \omega_{a}\\
\label{rs3}
\Bigl(\frac{F'}{\Mpc}\Bigr)'(z)
+ c\Bigl(\frac{F}{\Mpc}\Bigr)'(z) + 2 \lambda = 0. 
\end{gather}
Now \eqref{rs2} implies that for all $a$, $(\pm g_{a},\pm\omega_a)$ is
K\"ahler--Einstein and
\begin{equation} \label{rs4}
\frac{F'(z)}{\Mpc(z)} + c\frac{F(z)}{\Mpc(z)} = 2s_{a} - 2 \lambda 
(z + 1/x_{a}).
\end{equation}
Conversely this implies~\eqref{rs2}--\eqref{rs3}, the latter being just the
derivative of \eqref{rs4}.

As in \cite[\S2.4]{ACGT4}, since $\Theta(z)=F(z)/\Mpc(z)$, an application of
l'H\^opital's rule shows that~\eqref{boundary} is equivalent to
\begin{equation}\label{boundaryFbis}
F(\pm 1)=0,\qquad \Psi(-1)=2(d_0+1),\qquad
\Psi(1)=-2(d_\infty+1),
\end{equation}
where $F'(z)=\Psi(z)\Mpc(z)$.  Hence evaluating~\eqref{rs4} at $z=\pm1$, we
have
\begin{align}
\label{lambda}
2\lambda &= d_{0} + d_{\infty}+2\\
\label{kahlerclasses}
2 s_{a} x_{a} &= (d_{\infty}+1)(1-x_{a}) + (d_{0}+1)(1+x_{a}),
\end{align}
both expressions being manifestly positive (so the base manifolds $S_a$ have
positive scalar curvature). These equations allow us to rewrite~\eqref{rs4} as
a single equation
\begin{equation} \label{rs5}
\frac{F'(z)}{\Mpc(z)} + c\frac{F(z)}{\Mpc(z)} = (d_{0}+1)(1-z)
 -(d_{\infty}+1)(1+z)
\end{equation}
and~\eqref{kahlerclasses}--\eqref{rs5} imply~\eqref{rs4}. Using~\eqref{rs5},
the boundary conditions~\eqref{boundaryFbis} reduce to
\begin{equation} \label{rsboundary}
F(\pm 1) = 0.
\end{equation}
Hence we must solve~\eqref{kahlerclasses}--\eqref{rsboundary} subject to
$0<|x_a|<1$ and $F(z)>0$ for $z\in (-1,1)$. Clearly~\eqref{kahlerclasses}
gives $x_a=(d_{0} + d_{\infty}+2)/(2s_a+d_\infty-d_0)$ and so we must have
\begin{align}
\label{linebdlcond1}
s_{a} &> d_{0} + 1 && \text{if} \quad \omega_{a}>0, \\
\label{linebdlcond2}
s_{a} &< -(d_{\infty} +1) &&  \text{if} \quad \omega_{a} < 0.
\end{align}
Restricting the formula~\eqref{rho} for $\rho$ to the zero and infinity
sections $\zE_{0}$ and $\zE_{\infty}$, we see that these are actually necessary
conditions for $c_1(M)=[\rho/2\pi]$ to be positive.

We now observe that
\begin{equation} \label{Fsolution}
F(z) = e^{-cz}\int_{-1}^{z}e^{ct}\bigl((d_{0}+1)(1-t)-(d_{\infty}+1)(1+t)
\bigr)\Mpc(t) dt
\end{equation}
solves \eqref{rs5} and \eqref{rsboundary} iff $G(c)=0$, where
\begin{equation*}
G(k)= \int_{-1}^{1}e^{kt}((d_{0}+1)(1-t)- (d_{\infty}+1)(1+t))\Mpc(t) dt
= e^{kt_0}\int_{-1}^{1} e^{k(t-t_0)}(t-t_{0})g(t) dt
\end{equation*}
for some $t_{0}\in (-1,1)$ and $g(t)$ with $g<0$ on $(-1,1)$. Clearly
$e^{-kt_{0}} G(k)$ is a strictly decreasing function of $k$ tending to
$\mp\infty$ as $k\to\pm\infty$, so it has a unique zero $c$ (consistent with
the uniqueness of Ricci solitons). Since $F'$ has exactly one zero (namely
$t_{0}$) in $(-1,1)$, $F(\pm1)=0$ and $F$ is positive near the endpoints, it
is positive on $(-1,1)$. We deduce the following equivalence, essentially due
to Koiso~\cite{koiso}.

\begin{thm}\label{ricci-solitons} Let $S=\prod_{a\in\cA} S_a$ be a finite
product $(\cA\subset\Z^+)$ of compact K\"ahler--Einstein manifolds $(S_a,\pm
g_{a},\pm\omega_a)$ with scalar curvatures $\Scal_a = \pm 2d_a s_a$ and let
$M= P(\cO\otimes \C^{d_0+1} \oplus \cL\otimes \C^{d_{\infty}+1}) \to S$,
where $\cL=\bigotimes_{a\in\cA} \cL_a$ and $\cL_a$ are line bundles over
$S_a$ with $c_1(\cL_a)=[\omega_a/2\pi]$. Then the following conditions are
equivalent\textup:
\begin{itemize}
\item the conditions \eqref{linebdlcond1}--\eqref{linebdlcond2} are
satisfied\textup;
\item $(M,J)$ is a Fano manifold\textup; 
\item there exists a K\"ahler--Ricci soliton on $(M,J)$.
\end{itemize}
In this case, the K\"ahler--Ricci soliton $(g,\omega)$ is admissible with
$\lambda =(d_{0} + d_{\infty}+1)/2$ and $V= (c/2) {\grad}_g z$ for a suitable
real constant $c$.
\end{thm}
Our arguments and the fact that any Fano manifold is simply connected show
that Theorem~\ref{ricci-solitons} gives all compact K\"ahler--Ricci solitons
compatible with a hamiltonian 2-form of order 1 as above. We also have the
following standard corollary.
\begin{cor}\label{existence-KE}\cite{koiso}
Let $M^{2m}= P(\cO\otimes \C^{d_0+1} \oplus \cL\otimes \C^{d_{\infty}+1}) \to
S$, as in the above theorem.  Then there is a K\"ahler--Einstein metric on $M$
if and only if the conditions~\eqref{linebdlcond1}--\eqref{linebdlcond2} are
satisfied and the Futaki invariant $\mathfrak{F}_{[\rho]}(K)$ vanishes.
\end{cor}
The Futaki invariant $\mathfrak{F}_{[\rho]}(K)$ is a nonzero multiple of the
coefficient of $z^{m+2}$ in the extremal polynomial ${F}_{[\rho]}(z)$ as
defined in~\cite{ACGT4} (which is the leading coefficient if it is nonzero).
Hence its vanishing is equivalent to ${F}_{[\rho]}$ having degree at most
$m+1$. Unfortunately, verifying this condition is not easy (it leads to a
non-trivial diophantine problem); we will rediscover some K\"ahler--Einstein
examples of \cite{koi-sak1, koi-sak2} in the next section as a byproduct of
our study of WBF metrics.

\section{Constructions of WBF K\"ahler metrics}
\label{WBFsec}

We turn now to the construction of admissible WBF K\"ahler metrics on
admissible projective bundles.  By Proposition~\ref{equations}, an admissible
metric $g$ with $F(z)=\Theta(z)\Mpc(z)$, $\Mpc(z)=\prod_a (1+x_az)^{d_a}$
($0\leq a\leq \infty$, $d_a\geq 0$) is WBF, with $\tilde\rho$ a linear
combination of the hamiltonian $2$-form $\phi$ and the K\"ahler form $\omega$,
precisely when the metrics $g_a$ are K\"ahler--Einstein and
\begin{equation} \label{Fprime}
F'(z)=\Mpc(z)Q(z)
\end{equation}
for a polynomial $Q$ of degree $\leq 2$ with
\begin{equation} \label{einsteinconst}
Q(-1/x_a)=2s_a \qquad (a\in\smash{\hat\cA}).
\end{equation}
In this case $F$ is the extremal polynomial of the corresponding admissible
K\"ahler class \cite{ACGT4} and the WBF K\"ahler metric is K\"ahler--Einstein
iff $Q$ has degree $\leq 1$.

Since $g$ is, in particular, extremal, we know from \cite{ACGT4} (and it is
straightforward to check) that the positivity \eqref{positivity} and endpoint
conditions \eqref{boundary} may be replaced with
\begin{align}
\label{positivityF}
F &> 0 \quad\text{on}\quad (-1,1)\\
\label{boundaryF}
F(\pm 1) &= 0, \qquad F'(\pm 1) = \mp 2\Mpc(\pm1).
\end{align}

Using equations \eqref{Fprime} and \eqref{einsteinconst}, equation
\eqref{boundaryF} implies that $Q(-1)=2(d_0+1)$ and $Q(1)=-2(d_\infty+1)$. We
remark that since $Q(z)$ therefore changes sign only once on $(-1,1)$, so does
$F'(z)$ (since $\Mpc(z)$ is positive). Hence $F(z)$ (and $F(z)/\Mpc(z)$) will
be positive on $(-1,1)$ as soon as~\eqref{boundaryF} is satisfied.

The general quadratic $Q$ satisfying $Q(-1)=2(d_0+1)$ and
$Q(1)=-2(d_\infty+1)$ is
\begin{equation}
Q(z) = B(1-z^2)+(d_0+1)(1-z)-(d_\infty+1)(1+z)
\end{equation}
(and the K\"ahler--Einstein case is when $B=0$).
Equation~\eqref{einsteinconst} gives
\begin{equation*}
2s_a x_a^2 = B(x_a^2-1) + (d_0+1)(1+x_a)x_a+(d_\infty+1)(1-x_a)x_a.
\end{equation*}
We write $B=B_a$ for the solutions of these equations $(a\in\cA)$, so that
\begin{equation}\label{ca}
B_a:=x_a \bigl((d_0+1)(1+x_a)+(d_\infty+1)(1-x_a)-2 s_a x_a\bigr)/ (1-x_a^2).
\end{equation}

On the other hand, given the above, then \eqref{boundaryF} is satisfied iff we
set $F(z)=\int_{-1}^z \Mpc(t)\bigl(B (1-t^2) +(d_0+1)(1-t)-(d_\infty+1)(1+t)
\bigr) dt$ and
\begin{equation} \label{integralbc}
\int_{-1}^1 \Mpc(t)\bigl(B(1-t^2)+(d_0+1)(1-t)-(d_\infty+1)(1+t)\bigr)dt=0.
\end{equation}
Since $\Mpc(t)(1-t^2)$ is positive on $(-1,1)$, this determines $B$ uniquely,
once all other quantities are known.  Hence, in order to complete the
construction, we must show that $B=B_a$ solves~\eqref{integralbc} for all
$a\in\cA$.  Multiplying by $1-x_a^2$, this means that $h_a=0$ for all such
$a$, where
\begin{multline}\label{hadef}
h_a=\int_{-1}^1 \Mpc(t)\Bigl(
(1-x_a^2)\bigl((d_0+1)(1-t)-(d_\infty+1)(1+t)\bigr)\\
+x_a\bigl((d_0+1)(1+x_a)+(d_\infty+1)(1-x_a)-2s_ax_a\bigr)(1-t^2)\Bigr)dt.
\end{multline}
Our strategy for solving this problem is to use the equations
$\{h_a=0:a\in\cA\}$ to determine $\{x_a:a\in\cA\}$ as functions of
$\{s_a:a\in\cA\}$. For given $s_a=p_a/q_a$, we obtain a WBF K\"ahler metric on
the corresponding projective bundle iff we can find solutions $x_a$ with
$0<|x_a|<1$. We note that $h_a=\int_{-1}^1 \Mpc(t)k_a(t)dt$, where
\begin{equation}\label{kadef}
k_a(t)=\bigl((d_0+1)(1+x_a)(1-t)-(d_\infty+1)(1-x_a)(1+t)\bigr)(1+x_a t)
-2 s_ax_a^2(1-t^2).
\end{equation}
We remark that if $s_b\neq s_a$, $x_b$ cannot equal $x_a$, since $\int_{-1}^1
\Mpc(t)(1-t^2) dt$ is positive.  Hence if $x_a=x_b$, then $s_a=s_b$ and
$S_a\times S_b$ is K\"ahler--Einstein. Thus we do not need to check that $x_a$
are distinct: if $x_a=x_b$, we still get a WBF K\"ahler metric, but the
hamiltonian $2$-form has fewer constant roots.

Note also that we can replace the momentum coordinate $z$ by $-z$: this allows
us to replace $s_a$ by $-s_a$ and $x_a$ by $-x_a$, provided we interchange
$d_0$ and $d_\infty$.

\begin{rem} If the base manifolds are all $\C P^{d_a}$ and come in
pairs with equal dimensions with $d_0=d_\infty$ and (say) $d_{2k-1}=d_{2k}$
for $k\geq 1$, then it is straightforward to find some K\"ahler--Einstein
solutions to the equations $h_a=0$ by symmetry: for $|q_a|<(d_a+1)/(d+1)$ with
$q_{2j-1}=-q_{2j}$, set $s_a=(d_a+1)/q_a$ and $x_a=q_a (d+1)/(d_a+1)$; then
the integrand defining $h_a$ is an odd function of $t$, hence $h_a=0$. These
metrics are special cases of those of Koiso--Sakane~\cite{koi-sak1,koi-sak2}
and provide examples where the necessary and sufficient conditions of
Corollary~\ref{existence-KE} are verified (see also
Corollary~\ref{koiso-sakane-examples} below).
\end{rem}

\subsection{WBF K\"ahler metrics over a K\"ahler--Einstein manifold}

Let us consider the case when the base is a single K\"ahler--Einstein manifold
i.e., $\#\cA=1$. In the absence of blow-downs, this case was also considered
in \cite{ACGT3}. Dropping the $a$ subscript for this unique $a\in\cA$, we may
assume that we have to find $0<x<1$ such that $h(x)=0$, where
\begin{equation*}\begin{split}
h(x)&=
\int_{-1}^1 (1+t)^{d_0}(1-t)^{d_\infty}(1+x t)^{d} k(x,t)dt\\
k(x,t)&= \bigl((d_0+1)(1+x)(1-t)-(d_\infty+1)(1-x)(1+t)\bigr)(1+x t)
-2 sx^2(1-t^2).
\end{split}\end{equation*}
(Alternatively we could assume that e.g., $d_{0}\leq d_{\infty}$,
but then both $x$ positive and $x$ negative have to be considered.)  Since
$(1+t)^{d_0+1}(1-t)^{d_\infty+1}(1+xt)^{d+1}(1-xt)$ vanishes at $t=\pm1$ we
may add its derivative onto the integrand to obtain
\begin{equation}\label{second}\begin{split}
h(x)&=
\int_{-1}^1 (1+t)^{d_0+1}(1-t)^{d_\infty+1}(1+x t)^{d} x \hat k(x,t)dt\\
\hat k(x,t)&= (d_0+d_\infty+2-d)(1+xt)+2x((d+1)t-s).
\end{split}\end{equation}
Using the two integral formulae for $h(x)$, we make the following observations:
\begin{itemize}
\item $h(1)$ has sign $(d_0+1)-s$;
\item if $d\neq d_0+d_\infty+2$, $h(x)$ has sign $d_0+d_\infty+2-d$ for 
$x$ small and positive;
\item if $d=d_0+d_\infty+2$, then $h(x)$ has sign $(d+1)(d_0-d_\infty)-s(d+2)$
(if this is nonzero) for small nonzero $x$.
\end{itemize}
For this last case, evaluating $h(x)/x^2$ at $x=0$ gives $(s+(d+1))I_0+
(s-(d+1))I_\infty$ where $I_0$ and $I_\infty$ are integrals related by the
identity $(d_0+2)I_0=(d_\infty+2)I_\infty$.

If $d=d_0+d_\infty+2$ and $(d+1)(d_0-d_\infty)=s(d+2)$, it is easy to see
(integrating~\eqref{second} by parts) that there are no solutions of $h(x)=0$
with $0<x<1$.

Since $h$ is continuous, these sign observations lead to existence results.

\begin{thm} \label{t:wbf1}
Let $(S,g_S,\omega_S)$ be a compact Hodge K\"ahler--Einstein $2d$-manifold of
scalar curvature $2ds$ and let $E_{0}$, $E_{\infty}$ be projectively-flat
hermitian vector bundles of ranks $d_{0}+1$, $d_{\infty}+1$ over $S$ with with
$\cb_1(\vE_\infty)-\cb_1(\vE_0)=[\omega_S/2\pi]$.  Then there is an admissible
weakly Bochner-flat K\"ahler metric on $P(E_{0} \oplus E_{\infty})\to S$
when\textup:
\begin{itemize}
\item $S$ has nonpositive scalar curvature $(s\leq 0)$, $d\geq
d_0+d_\infty+2$, unless $d=d_0+d_\infty+2$ and $(d+1)(d_\infty-d_0)\leq
|s|(d+2)$\textup;

\item $S$ has positive scalar curvature $(s>0)$, $(d_0+1) > s$, and $d\geq
d_0+d_\infty+2$, unless $d=d_{0}+d_{\infty} +2$ and $d_{0} >
d_{\infty}$\textup;

\item $S$ has positive scalar curvature $(s>0)$, $(d_0+1)<s$, and $d<
d_0+d_\infty+2$.

\end{itemize}
\end{thm}

When $d_0=d_\infty=0$ and $S$ is a positive K\"ahler--Einstein manifold, these
existence results are sharp. In particular, when $S=\C P^d$, we obtain the
following result.

\begin{thm}\label{WBFunique}
There is a weakly Bochner-flat K\"ahler metric on $P(\cO \oplus \cO(q) ) \to
\C P^d$ with $q>0$ if and only if $d=1$ and $q=1$ or $d \geq 2$ and
$q>d+1$. The weakly Bochner-flat K\"ahler metric is then unique up to
automorphism and scale.
\end{thm}
\begin{proof} Any WBF K\"ahler metric is extremal and the extremal K\"ahler
metrics on $M=P(\cO \oplus \cO(q))\to \C P^d$ have cohomogeneity one under a
maximal compact connected subgroup of ${\Aut}(M,J)$~\cite{calabi1}. Since any
two such subgroups are conjugate in the connected component $\Aut(M,J)^0$, it
follows that, up to pullback by a automorphism, the WBF K\"ahler metrics on
these manifolds must be admissible.  The existence of a WBF K\"ahler metric in
the stated cases follows from Theorem~\ref{t:wbf1} above, so it remains to
establish the nonexistence and uniqueness results.

For the case $d=1$, we compute that
\begin{equation}\label{m1h}
h(x)=\tfrac{4}{3} x\bigl( x^2+1 - 2sx \bigr)
\end{equation}
and clearly there is a unique solution $0<x<1$ to $h(x)=0$ iff $s>1$. Since
$S$ in this case is $\C P^1$, $\cK^{-1}=\cO(2)$ and the only possibility is
$s=2$, $\cL=\cO(1)$, in accordance with the classification of~\cite{ACG1}.

For the case $d=2$ we calculate directly that
\begin{equation}
\label{neededfor6mnfclass}
h(x)=\tfrac{8}{15}x^2 \bigl( 6x - s (x^2+5) \bigr)
\end{equation}
and clearly there is a unique solution $0<x<1$ to $h(x)=0$ iff $0<s<1$.

We now assume $d\geq 3$ and compute the integral (e.g., by substitution) to
get:
\begin{multline*}
-\tfrac12(d+1)(d+2)(d+3)x^2 \,h(x) =\\
(1-x)^{d+2}
\bigl(d+1 + ((d+1)(d+2) + 2 s) x + ((d+1)(d+3) + 2(d+2)s )x^2\bigr)\\
- (1+x)^{d+2}
\bigl(d+1 - ((d+1)(d+2) - 2 s) x + ((d+1)(d+3) - 2(d+2)s )x^2\bigr).
\end{multline*}
If $x=(y-1)/(y+1)$ and $f(y)= -(d+1)(d+2)(d+3)(y+1)^{d+1}(y-1)
h(x)/2^{d+4}$ then
\begin{multline*}
f(y)= (d+1)(s+1) - (d+2)(d + 1 + 2s) y + (d+3)(d + 1 + s) y^2\\
+ y^{d+2}
\bigl( -(d+3) (d + 1 - s) + (d+2)(1 + d - 2 s) y + (d+1)(s-1) y^2 \bigr).
\end{multline*}
The zeros of $h(x)$ in $(0,1)$ correspond to the zeros of $f(y)$ in
$(1,\infty)$.  The latter problem is more amenable to calculus, since
$f(1)=f'(1)=f''(1)=0$ and $f'''(y)= (d+1)(d+2)(d+3) y^{d-1} P(y)$, where
\begin{equation*}
P(y)= -d (1 + d - s) + (d+2)(d+1 - 2 s) y + (d+4)(s-1) y^2.
\end{equation*}
Now $P(1)=d-2$, which is positive for $d>2$, while $P(0)$ is nonpositive since
$s \leq d + 1$. Hence $P(y)$ is positive in $(1,\infty)$ unless $s<1$, in
which case it has a unique zero. If $P(y)$ is positive in $(1,\infty)$, then
so is $f'''$, hence $f''$, $f'$ and $f$, because we know that
$f(1)=f'(1)=f''(1)=0$. This gives the nonexistence.  Similarly, when $f'''(y)$
has a unique zero in $(1,\infty)$, so does $f$, which gives the required
uniqueness.
\end{proof}

Note that the proof above in the case $d=2$ also gives us the following
result.

\begin{thm}\label{wbflinebdlover4mnf} Let $S$ be a compact K\"ahler--Einstein
complex surface.  There is an admissible weakly Bochner-flat K\"ahler metric
with $\#\cA=1$ on $P(\cO \oplus \cL) \rightarrow S$ if and only if $S$ is a
positive K\"ahler--Einstein manifold and $\cL = \cK^{-q/p}$, where integers
with $|q|>p>0$ such that $\cK^{-1/p}$ is the primitive ample root of the
canonical bundle of $S$. The admissible weakly Bochner-flat K\"ahler metrics
is then unique up to automorphism and scale.
\end{thm}

We end this paragraph by studying in more detail the case $d=1$ and
$d_{0}+d_\infty=1$, when $M$ is a $\C P^2$-bundle over a compact Riemann 
surface $S_{1} = \Sigma$.  Again, 
we assume
without loss that $0<x<1$.

When $d_{0}=1$ and $d_{\infty}=0$ we have $h(x) = 0$ iff $r(x)= (3 - s) x^{2}
+ (4 - 5s) x + 5 = 0$.  If $r(x)=0$ then $s \geq 4/5$ and the (positive
definite) metric $g_{\Sigma}$ is a constant curvature metric on $\Sigma=\C
P^{1}$, so we must have that $E_{0}= \cL_{0} \otimes \C^{2}$, $E_{\infty} =
\cL_{\infty} \otimes \C$ for some line bundles $\cL_{0}$, $\cL_{\infty}$ and
that $\omega_{\Sigma}/2\pi$ is integral.  Thus $s=1$ or $s=2$ (since $s= 2/q$
for $q \in \Z^{+}$). However $r(x)$ does not have a root in $(0,1)$ in either
case.

When $d_{0}=0$ and $d_{\infty}=1$, we have $h(x) = 0$ iff $r(x)= (3 + s) x^{2}
- (4 + 5s) x + 5 = 0$. Since $r(x)= (1-x) (5-4x) > 0$ for $s=1$ and
$\frac{\partial}{\partial s}r(x) = x(x-5) < 0$, $r(x)$ has no roots in $(0,1)$
for $s \leq 1$.  Then we may assume that $S_{1}=\C P^{1}$ and that $E_{0} =
\cO \otimes \C = \cO$ and $E_{\infty} = \cL \otimes \C^{2}$, where $\cL$ is a
holomorphic line bundle with $c_{1}(\cL) = [\omega_{\Sigma}/2\pi]$. By
integrality, the only possibility with $s>1$ is $s=2$, for which we find a
unique solution $x= (7-2\sqrt{6})/5$ in $(0,1)$ (so $\Sigma=\C P^{1}$ and
${\cL} = \cO (1)$). Observe that $B\neq 0$, so the corresponding metric is not
K\"ahler--Einstein.

\begin{thm}\label{wbfblowdown}
Let $E_{0}$, $E_{\infty}$ be projectively-flat hermitian vector bundles over a
compact Riemann surface $\Sigma$ with ranks $d_{0} + 1$, $d_{\infty}+ 1$
respectively.  Then there is an admissible weakly Bochner-flat K\"ahler metric
on $P(E_{0} \oplus E_{\infty}) \to \Sigma$ with $d_{0}+d_{\infty}=1$ if and
only if \textup(without loss\textup) $d_{0}=0$, $d_{\infty}=1$, $\Sigma=\C
P^{1}$ and $E_{0} = \cO$ while $E_{\infty} = \cO(1) \otimes \C^{2}$.  The
admissible weakly Bochner-flat K\"ahler metric is then unique up to
automorphism and scale.
\end{thm}

\subsection{WBF K\"ahler metrics over a product of K\"ahler--Einstein
manifolds}

In this paragraph and the next, we consider the case that $d_{0}=d_{\infty}=0$
and $\#\cA=2$ in detail. We will assume that the base $S$ is a global product
of two K\"ahler Einstein manifolds $S_a$ ($a=1,2$) of dimensions $2d_a>0$. We
postpone a detailed discussion of the case $d_{1}=d_{2}=1$ to the next
paragraph (where we also consider the case where $S$ is a local product). In
this setting we have (up to a constant factor)
\begin{align*}\notag
h_1(x_1, x_2)&=
\int_{-1}^1 (1+x_1 z)^{d_1}(1+x_2 z)^{d_2}
              \bigl( x_1 (x_1 s_1 - 1) (1-z^2) + z (1-x_1^2) \bigr) dz,\\
h_2(x_1, x_2)&=
\int_{-1}^1 (1+x_1 z)^{d_1}(1+x_2 z)^{d_2}
              \bigl( x_2 (x_2 s_2 - 1) (1-z^2) + z (1-x_2^2) \bigr) dz.
\end{align*}
We are looking for common zeros of these functions with $0<|x_a|<1$.  Let us
note what we know about these functions on the boundary of this domain:
\begin{itemize}
\item when $x_1=0$, $h_1$ has the same sign as $x_2$;
\item when $x_1=\pm 1$, $h_1$ has the same sign as $s_1\mp 1$;
\item when $x_2=0$, $h_2$ has the same sign as $x_1$;
\item when $x_2=\pm 1$, $h_2$ has the same sign as $s_2\mp 1$.
\end{itemize}
In particular, the curves $h_1=0$ and $h_2=0$ both pass through $(0,0)$ and we
know the gradients of these curves at $(0,0)$, since $\partial h_a/\partial
x_a = 2(d_a-2)/3$ and $\partial h_a/\partial x_b = 2d_b/3$ for $b\neq
a$. Hence along $h_1=0$ we have $dx_2/dx_1=(2-d_1)/d_2\leq 1$ at $(0,0)$,
while along $h_2=0$ we have $dx_1/dx_2=(2-d_2)/d_1\leq 1$ at $(0,0)$ so that
$dx_2/dx_1 = d_1/(2-d_2)$ (infinite when $d_2=2$).  Furthermore, if both
curves have negative gradients, $dx_2/dx_1$, at $(0,0)$---that is, both curves
emanate from the origin into the fourth quadrant---then we must have that
$d_1>2$ and $d_2>2$. Hence the difference in the gradients, namely
$2(d_1+d_2-2)/d_2(d_2-2)$, is positive, so that the curve $h_1=0$ is above the
curve $h_2=0$ for $x_1>0$ near $(0,0)$.

There are two separate types of solutions to seek: those with $x_1$ and $x_2$
of opposite sign, and those with $x_1$ and $x_2$ of the same sign.  Figures
1--2 plot examples of the graphs of $h_{1}=0$ (solid) and $h_{2}=0$ (dashed)
in each case.

We consider first the case of opposite signs, and without loss, we seek
solutions with $x_1>0$ and $x_2<0$. Suppose now that $s_1>1$ and
$s_2<-1$. Then
\begin{itemize}
\item $h_1$ changes sign on any path from $x_1=0, x_2<0$ to $x_1=1, x_2\leq
0$;
\item $h_2$ changes sign on any path from $x_2=0, x_1>0$ to $x_2=-1, x_1\geq
0$.
\end{itemize}
It follows by continuity that the curves $h_1=0$ and $h_2=0$ must cross.
\begin{lemma} \label{indef}
If $s_1>1$ and $s_2<-1$ then there exist $x_1\in (0,1)$, $x_2\in(-1,0)$ such
that $h_1(x_1,x_2)=0=h_2(x_1,x_2)$.
\end{lemma}
\begin{proof} Since $h_1$ is negative on the half-line $(x_1=0,x_2<0)$
and positive on $x_1=1$, there is a connected component ${\cC}$ of the
curve $h_1=0$ in the square $(0,1)\times [0,-1]$ which crosses $x_2=-1$ for
some $x_1\in (0,1)$, and it either crosses $x_2=0$ for some $x_1\in (0,1)$, or
it emanates from the origin, and, within the square, is initially above the
curve $h_2=0$, as in Figure 1. It follows that $h_2$ changes sign
on ${\cC}$, hence vanishes by continuity and connectedness.
\end{proof}

\centerline{\includegraphics[height=1in,width=1in]{233neg2.epsf}
\qquad\qquad\qquad \includegraphics[height=1in,width=1in]{122neg2.epsf}}
\vspace{0.2cm}
\par\centerline{Figure~1:~$d_{1}=2$, $d_{2}=3$, $s_{1}=3$,
$s_{2}=-2$ and $d_{1}=1$, $d_{2}=2$, $s_{1}=2$,
$s_{2}=-2$}

\vspace{0.4cm}

Let us turn now to the case that $x_1$ and $x_2$ have the same sign, so
without loss, $x_1>0$ and $x_2>0$.  Suppose that $s_1<1$ and $s_2<1$. Then
\begin{itemize}
\item $h_1$ changes sign on any path from $x_1=0, x_2>0$ to $x_1=1, x_2\geq
0$;
\item $h_2$ changes sign on any path from $x_2=0, x_1>0$ to $x_2=1, x_1\geq
0$;
\item the curve $h_1=0$ lies below the line $x_1=x_2$ for $x_1>0$ near
$(0,0)$, and is strictly below unless $d_1=1$;
\item the curve $h_2=0$ lies above the line $x_1=x_2$ for $x_2>0$ near
$(0,0)$, and is strictly above unless $d_2=1$.
\end{itemize}
Again we see that the curves $h_1=0$ and $h_2=0$ must cross, except perhaps in
the case $d_1=d_2=1$, which we shall consider in the next paragraph.

\begin{lemma} \label{2ndderiv}
If $s_1<1$ and $s_2<1$, and $d_1,d_2$ are not both $1$, then there exist
$x_1,x_2\in (0,1)$ such that $h_1(x_1,x_2)=0=h_2(x_1,x_2)$.
\end{lemma}
\begin{proof} As in the previous lemma, there is a connected component
${\cC}$ of the curve $h_1=0$ in the square $(0,1)\times [0,1]$ which
crosses $x_2=1$ for some $x_1\in (0,1)$, and it either crosses $x_2=0$ for
some $x_1\in (0,1)$, or it emanates from the origin. In the latter case, we
need to know that $h_1=0$ is initially below $h_2=0$, so that $h_2$ is
initially positive. Since not both $d_{1}$ and $d_{2}$ equal one, this follows
from the observations prior to the statement of the lemma.
\end{proof}

\centerline{\includegraphics[height=1in,width=1in]{132over34over5.epsf}}
\par\centerline{Figure~2:~$d_{1}=1$, $d_{2}=3$, $s_{1}=2/3$, $s_{2}=4/5$}
\vspace{.2cm}

Let us summarize what we have established, excluding the case $d_1=d_2=1$.

\begin{thm}  
Let $S_a$ $(a=1,2)$ be compact K\"ahler--Einstein $2d_a$-manifolds with
$d_a\geq 1$ not both one. Let $\cK_a$ be the canonical bundles, and suppose
\textup(without loss unless $S_a$ is Ricci-flat\textup) that the K\"ahler form
$\pm\omega_a$ is integral.  Let $\cL_a$ be line bundles on $S_a$ with
$c_1(\cL_a)=[\omega_a/2\pi]$ and, if $S_a$ is not Ricci-flat, let $\cL_a$ be
$\cK_a^{-q_a/p_a}$ tensored by a flat line bundle, for integers $p_a,q_a$
where $\cK^{-1/p_a}$ is the primitive ample root of the canonical bundle of
$S_{a}$.  Then there is an admissible weakly Bochner-flat K\"ahler metric on
$P(\cO\oplus\cL_1\otimes\cL_2)\to S_1\times S_2$ in the following
cases\textup:
\begin{itemize}
\item $S_1$ and $S_2$ have positive scalar curvature, $0<q_1<p_1$ and
$0<-q_2<p_2$\textup;
\item for $a=1,2$, $q_a>p_a$ if $S_a$ has positive scalar curvature, 
$q_a>0$ if $S_a$ has negative scalar curvature, and $\omega_{a}$ is 
positive  if $S_{a}$ is Ricci flat\textup.
\end{itemize}
\end{thm}

\begin{cor} There is a weakly Bochner-flat K\"ahler metric on
$P(\cO\oplus\cO(q_1,q_2))\to\C P^{d_1}\times\C P^{d_2}$ in the following
cases\textup:
\begin{itemize}
\item $q_1>d_1+1$ and $q_2>d_2+1$\textup;
\item $1\leq q_1\leq d_1$ and $1\leq -q_2\leq d_2$.
\end{itemize}
\end{cor}
We will see in the next paragraph that this corollary also holds for
$d_{1}=d_{2}=1$.  We conjecture that all WBF K\"ahler metrics on
$P(\cO\oplus\cO(k_1,k_2))\to\C P^{d_1}\times\C P^{d_2}$ are given by this
corollary and that the metric is unique (up to automorphism and scale) in
each case. As in Theorem \ref{WBFunique}, extremal K\"ahler metrics on these
manifolds are cohomogeneity one, hence of linear type, but unless
$d_{1}=d_{2}=1$ (see next paragraph) we have not been able to establish the
relevant nonexistence and uniqueness results for solutions of $h_1=0=h_2$.

We note also that if $d_1=d_2$ (including the case $d_1=d_2=1$) and $k_1=-k_2$
in the above corollary, we have not just a WBF K\"ahler metric, but a
K\"ahler--Einstein metric, as found by Koiso and
Sakane~\cite{koi-sak1,koi-sak2}.

\begin{cor} \label{koiso-sakane-examples} \cite{koi-sak1,koi-sak2}
On $P(\cO\oplus \cO(q,-q)) \to \C P^d \times \C P^{d}$, with $1\leq q\leq d$,
there is a K\"ahler--Einstein metric, given \textup(on a dense open
set\textup) by
\begin{equation*}
g=\Bigl(\frac{d+1}{q}+z\Bigr)g_1 + \Bigl(\frac{d+1}{q}-z\Bigr)g_2
+ \frac{z^{2}-(d+1)^2/q^2}{F(z)} \,dz^{2} +
\frac{F(z)}{z^{2}-(d+1)^2/q^2} \,\theta^{2},
\end{equation*}
where $(g_1,\omega_1)$ and $(g_2,\omega_2)$ are Fubini--Study metrics on the
$\C P^d$ factors with holomorphic sectional curvature $2/q$, $d\theta=
\omega_1-\omega_2$ and $F(z)= \int_{-1}^{z}
2t\bigl(\frac{(d+1)^2}{q^2}-t^{2}\bigr)\,dt =
-\smash{\frac{(d+1)^2}{q^2}}(1-z^2) + \frac12 (1-z^4)$.
\end{cor}
\begin{proof}
Let $s_1=-s_2=\frac{d+1}{q}$ and $x_1=-x_2=\frac{q}{d+1}$.  Then clearly
$h_1(x_1,x_2)=h_2(x_1,x_2)=0$. Further, $x_a=1/s_a$ so the WBF metric is
K\"ahler--Einstein.
\end{proof}

\subsection{WBF K\"ahler metrics over a ruled surface}
\label{wbfoverruledsurfsec}

Let us now consider the case $d_1=d_2=1$, when the base is a product of
Riemann surfaces.  Thus we have
\begin{align*}\notag
h_1(x_1, x_2)&=
\int_{-1}^1 (x_1 z+1)(x_2 z+1)
              \bigl( x_1 (x_1 s_1 - 1) (1-z^2) + z (1-x_1^2) \bigr) dz\\
h_2(x_1, x_2)&=
\int_{-1}^1 (x_1 z+1)(x_2 z+1)
              \bigl( x_2 (x_2 s_2 - 1) (1-z^2) + z (1-x_2^2) \bigr) dz
\end{align*}
(up to a constant factor), which by integration gives
\begin{align*}
h_1(x_1,x_2)&=\tfrac{2}{15}
(5x_2-5x_1+10s_1x_1^{2}-7x_1^{2}x_2-5x_1^{3}+2s_1x_1^{3}x_2)\\
h_2(x_1,x_2)&=\tfrac{2}{15}
(5x_1-5x_2+10s_2x_2^{2}-7x_2^{2}x_1-5x_2^{3}+2s_2x_2^{3}x_1).
\end{align*}

Without loss, we look for solutions to $h_1(x_1,x_2)=0=h_2(x_1,x_2)$ with
$x_1>x_2$ and $x_1>0$.  Solving $h_1=0$ for $s_1$, we find that $s_1$ must be
positive, hence $s_1=2/q_{1}$ for some integer $q_{1}\geq 1$.  We then
establish the following three lemmas, the proofs of which can be found in
Appendix~\ref{appC}.

\begin{lemma} \label{s1is2} If $s_1=2$ then there exist
$(x_1,x_2) \in (0,1) \times (-1,1)$ such that $h_1(x_1,x_2)=h_2(x_1,x_2)=0$
iff $s_{2} \leq -2$.  Moreover, in this case the solution is unique. If $s_{2}
< -2$ the solution is in $(0,1) \times (0,1)$, i.e., $x_2>0$, while if
$s_{2}=-2$, the solution is $(\tfrac12,-\tfrac12)$.
\end{lemma}

\begin{lemma} \label{s1is1} If $s_1=1$ then there exist
$(x_1,x_2) \in (0,1) \times (-1,1)$ such that $h_1(x_1,x_2)=h_2(x_1,x_2)=0$
iff $s_{2}<-1$.  Moreover, in this case the solution is unique and
$x_{2} > 0$.
\end{lemma}

\begin{lemma} \label{s1is2overk} If $s_1=2/q_{1}$, where $q_{1}\in\Z$ and 
$q_{1}\geq 3$, then there exist $(x_1,x_2) \in (0,1) \times (-1,1)$ such that
$h_1(x_1,x_2)=h_2(x_1,x_2)=0$ \textup(with $s_2=2/q_2$ if $x_2>0$\textup) iff
$-s_{1}< s_{2}<1$.  Moreover, in this case the solution is unique and
$x_{2} > 0$.
\end{lemma}

We do not need to assume $S_1$ and $S_2$ are compact for these arguments.
However, if $S_1$ is complete, it must be $\C P^{1}$ and the product $S_{1}
\times S_{2}$ is a (trivial) ruled surface. More generally we can suppose this
is the universal cover of compact K\"ahler surface, which is then a
geometrically ruled surface $S=P(E)$ over a Riemann surface $\Sigma$
with universal cover $S_2$. It is well-known that the existence of a local
product metric on $S$ is equivalent to $P(E)\to\Sigma$ admitting a flat
projective unitary connection. This in turn, by a famous result of Narasimhan
and Seshadri~\cite{ns}, is equivalent to polystability of $E$.  The above
lemmas therefore imply the following result.

\begin{thm}\label{wbfoverruledsurf}
Let $S$ be a Hodge $4$-manifold whose universal cover is a product of
constant curvature Riemann surfaces and suppose that $M=P(\cO\oplus\cL)\to
S$ has an admissible weakly Bochner-flat K\"ahler metric. Then
$S$ is a geometrically ruled surface $P(E)$ such that $E \to \Sigma$ is
polystable.  Let ${\mathbf f}, {\mathbf v} = c_1(VP(E)) \in
H^{2}(S,\Z)$ denote the classes of a fibre of $P(E) \to \Sigma$ and
of the vertical line bundle.  We then have $c_{1}(\cL) = (q_{1}/2) {\mathbf v}
+ q_{2}{\mathbf f}$ where $q_{1} \in \Z$, and $q_{2} \in \Z$ unless $q_{1}$ is
odd and $E \to \Sigma$ is not spin \textup(which may only happen when $\Sigma$
has genus $\mathbf g>1$\textup), in which case $q_{2} + 1/2 \in
\Z$. Furthermore, up to replacing $\cL$ by $\cL^{-1}$\textup:
\begin{itemize}  
\item if $\Sigma= \C P^{1}$, $S = \C P^{1}\times\C P^{1}$, and we
either have $q_{1}=1$ and $q_{2}=-1$, or we have $q_{1},q_{2} > 2$\textup;

\item if $\Sigma = T^{2}$, $q_{1}> 2$ and $q_{2}> 0$\textup;

\item if $\Sigma$ has genus $\mathbf g>1$, we either have $q_1 > 2$ and $q_2 >
q_1 ({\mathbf g} -1)$, or we have $q_{1}\in\{1,2\}$ and $0<q_{2}< q_1
({\mathbf g} -1)$.
\end{itemize}
Conversely, in each case there is a unique admissible weakly Bochner-flat
K\"ahler metric on $M$ up to automorphism and scale.
\end{thm}
Note that $E$ spin means that $\deg E$ is even. Since $\deg (E\otimes\cL)=\deg
E+2\deg\cL$, this condition (like polystability) is independent
of the choice of $E$ with $S=P(E)$.
\begin{proof}
We have seen already that $S=P(E)$ for $E\to\Sigma$ polystable.  If
$\Sigma = \C P^{1}$, $E$ is trivial and $S= \C P^{1} \times \C
P^{1}$.  If $\Sigma=T^2$, without loss $E$ is either $\cO \oplus \cL \to
\Sigma$ with $\deg \cL =0$ or the nontrivial extension of $\cO \to \Sigma$
\cite{suwa}. In either case $\deg E = 0$. Thus the non-spin case may only
happen when the genus of $\Sigma$ is at least $2$.

Let $\omega_{\C P^{1}}$ be the K\"ahler form of the Fubini--Study metric on
$\C P^{1}$ with volume one and let $\omega_{\Sigma}$ be a K\"ahler form of a
CSC K\"ahler metric on $\Sigma$ of volume one.

Let $\C P^{1} \times \tilde{\Sigma} \to S$ denote the universal cover
of $S$ (so $\tilde{\Sigma}$ covers $\Sigma$) and let $\pi_1\colon
\C P^{1} \times \tilde{\Sigma}\to \C P^{1}$ denote the projection to the first
factor. Then $\pi_1^* \omega_{\C P^{1}}$ descends to a closed $(1,1)$-form
on $S$ which represents ${\mathbf v}/2$, whereas ${\mathbf f} =
[\pi^* \omega_{\Sigma}]$.  Hence we see that a local product $q_{1} \omega_{\C
P^{1}} + q_{2} \omega_{\Sigma}$ corresponds to a line bundle $\cL \to
S$ with Chern class $(q_{1}/2) {\mathbf v} + q_{2}{\mathbf f}\in
H^{2}(S,\Z)$.  Now we note that $H^{2}(S,\Z) = \Z {\mathbf
h} \oplus \Z {\mathbf f}$, where ${\mathbf h} \in H^{2}(S,\Z)$
denotes the class of the dual of the ($E$-dependent) tautological line bundle
on $S$ (see e.g.,~\cite{fuj}).  Since ${\mathbf v} = 2 {\mathbf h} +
(\deg E) {\mathbf f}$, the integrality condition on $q_1,q_2$ for the
existence of $\cL$ follows immediately. Now we apply
Lemmas~\ref{s1is2}--\ref{s1is2overk}, bearing in mind that $s_{1}= 2/q_{1}$ and
$s_{2} = 2(1-{\mathbf g})/q_{2}$.
\end{proof}

\begin{cor}\label{P1timesP1} There is a weakly Bochner-flat K\"ahler metric
\textup(unique up to automorphism and scale\textup) on
$P(\cO\oplus\cO(q_1,q_2))\to\C P^{1}\times\C P^{1}$ if and only if $q_1>2$ and
$q_2>2$, or $q_1 = 1$ and $q_2 = -1$, the latter metric being
K\"ahler--Einstein.
\end{cor}
\begin{proof}
A WBF K\"ahler metric is in particular extremal and since extremal K\"ahler
metrics on these manifolds are cohomogeneity one, hence admissible (up to
automorphism), cf.~\cite{calabi1}, this follows from the above theorem and
Corollary~\ref{koiso-sakane-examples}.
\end{proof}

\subsection{WBF versus extremal K\"ahler metrics}

Any WBF K\"ahler metric is extremal, so our results provide examples of
extremal K\"ahler metrics in admissible K\"ahler classes in the sense
of~\cite{ACGT4}. By the results of~\cite{ACGT4}, we then obtain
$N$-dimensional families of such metrics near a WBF metric, where $N$ is the
number of K\"ahler--Einstein factors in the base. (In fact we do not need the
base metrics $g_a$ to be K\"ahler--Einstein to get an extremal K\"ahler
metric: it suffices in the above calculations that they are CSC and Hodge.)

\section{Classification of WBF K\"ahler metrics on compact $6$-manifolds}
\label{6-mnf-WBF}

Using the theory of~\cite{ACGT3,ACGT4}, the results of the previous
section yield the following classification result for compact
$6$-manifolds admitting WBF K\"ahler metrics.

\begin{thm}
Suppose that $(M,J,g,\omega)$ is a compact connected weakly Bochner-flat
K\"ahler $6$-manifold of order $\ell$. Then $\ell\in\{0,1\}$.

\begin{numlist}
  
\item If $\ell=0$ then $(M,J,g,\omega)$ is a local product of
K\"ahler--Einstein manifolds.

\item If $\ell=1$, then $(M,J)$ is biholomorphic to one of the following.

$(a)$ ${P}(\cO \otimes \cL) \rightarrow S$ where $S$ is a positive
K\"ahler--Einstein complex surface, $\cL = \cK^{-q/p}$ and $q>p>0$ are
integers and $p$ is the Fano index of $S$.

$(b)$ ${P}(\cO\oplus\cL)\to S$ where $S = P(E) \to \Sigma$
is a geometrically ruled surface such that $E \to \Sigma$ is polystable and
$\cL$ is given by Theorem~\textup{\ref{wbfoverruledsurf}}, excluding $P(\cO
\oplus \cO (-1,1)) \to \C P^1\times \C P^1$ \textup(which arises in
the case $\ell=0$ as it admits a K\"ahler--Einstein metric\textup).

$(c)$ $P(\cO \oplus \cO(1)\otimes\C^2)\to\C P^1$ \textup(a blow-down of
${P}(\cO \oplus \cO (1,-1))\rightarrow {\C P}^{1} \times {\C P}^{1}$\textup)

\smallbreak\noindent On each manifold in $(a)$--$(c)$, there is a unique
weakly Bochner-flat K\"ahler metric of up to automorphism and scale
\textup(and it has order $1$\textup).
\end{numlist}
\end{thm}

\begin{proof} In~\cite[Thm.~11]{ACGT4} we proved that a compact extremal
K\"ahler $6$-manifold admitting a hamiltonian $2$-form of order $2$ with the
extremal vector field tangent to the $\T^c$-orbits is isometric to $\C P^{3}$
with a Fubini--Study metric.  On the other hand, a compact K\"ahler
$6$-manifold with a hamiltonian $2$-form of order $3$ is biholomorphic to $\C
P^3$~\cite{ACGT3}, and hence, if it is extremal, it is again isometric to a
Fubini--Study metric.  Thus there are no compact WBF K\"ahler $6$-manifolds of
order $2$ or $3$.

Part (i) is immediate and the existence and biholomorphic classification in
part (ii) follow from Theorems~\ref{wbflinebdlover4mnf}, \ref{wbfblowdown} and
\ref{wbfoverruledsurf}. It remains to prove the uniqueness claim in (ii).  By
Theorems~\ref{wbflinebdlover4mnf}, \ref{wbfblowdown} and
\ref{wbfoverruledsurf}, and the well-known uniqueness result of Bando--Mabuchi
for K\"ahler--Einstein metrics~\cite{Bando-Mabuchi}, it suffices to prove that
any WBF K\"ahler metric is admissible (with the given bundle structures) up to
scale and automorphism, for which, using~\cite{ACGT3} again, it is enough to
show that the metric can be pulled back by an automorphism of $(M,J)$ so that
the extremal vector field $J \grad_{\smash g} \Scal_g$ becomes a nonzero
multiple of the generator of the canonical $S^1$-action. We now establish the
uniqueness in each case.

(a) By the classification of~\cite{tian}, $S$ is biholomorphic to $\C P^2$,
$\C P^1 \times \C P^1$, or a blow-up of $\C P^2$ at $k$ points in general
position for $3\leq k \leq 8$.  When $S= \C P^2$ or $S= \C P^1 \times \C P^1$,
the uniqueness follows from Theorem~\ref{WBFunique} and
Corollary~\ref{P1timesP1}, so it remains to consider the case that $S$ is a
blow-up of $\C P^2$. This has Fano index $p=1$, so $\cL = \cK^{-q}$ for
$q>1$. By Riemann--Roch, $H^0(S,\cL) \neq 0$ while $H^0(S, \cL^{-1})=0$ since
$\cL$ is not trivial.  Therefore, \cite[Props.~3--4]{ACGT4}
show that $M$ does not admit any CSC K\"ahler metrics.  In particular, any
other WBF K\"ahler metric $g'$ on $M$ must have order $1$ and is
therefore~\cite{ACGT3} admissible with respect to some ruling of $M$ over a
K\"ahler--Einstein surface $S'$ with $b_2(S')=b_2(M)-1=b_2(S)=4$. Since $g$
and $g'$ are both extremal, by~\cite{calabi1} we can assume, after pulling
back $g'$ by an automorphism, that ${\mathfrak i}_0(M,g')={\mathfrak
i}_0(M,g)$ in $\mathfrak h_0(M)$. Let $K, K'$ be the extremal vector fields of
$g$, $g'$.  Then $\cL_K\Scal_{g'}=\cL_{K'}\Scal_g=0$, so $K$ and $K'$ induce
hamiltonian Killing vector fields $X$, $X'$ on $S$, $S'$.  If either is zero,
$K'$ is a multiple of $K$ and we are done.  Otherwise,
$\mathfrak{h}(S),\mathfrak{h}(S) \neq 0$ so $S$, $S'$ are both (isomorphic to)
the blow-up of $\C P^2$ at three points.  The corresponding K\"ahler--Einstein
metrics agree up to automorphism and scale by~\cite{Bando-Mabuchi}, hence so
do $g$ and $g'$ (by Theorem~\ref{wbflinebdlover4mnf}).

(b)--(c) Here any K\"ahler class on $M$ is admissible, so $M$ admits no CSC
K\"ahler metrics by \cite[Thm.~5, Thm.~8 \& Rem.~8]{ACGT4} (in case (b), $\cL$
is, without loss, ample by Theorem~\ref{wbfoverruledsurf}).  Thus any WBF
metric on $M$ has order $1$. Being in an admissible class, its extremal vector
field must be a multiple of $K$ by \cite[Prop.~6]{ACGT4}.
\end{proof}

\begin{rem}
In the classification of WBF K\"ahler $4$-manifolds obtained in~\cite{ACG1}
the normalized Ricci form also has order 0 or 1. A naive dimension counting
argument~\cite{ACGT3} supports the conjecture that this feature persists in
higher dimensions. We also note that the base manifolds $S$ have Kodaira
dimension $-\infty$. In view of the examples of Theorem~\ref{t:wbf1}, this is
no longer true in dimension $\geq 8$.
\end{rem}

\appendix

\section{Proofs of Lemmas \ref{s1is2}, \ref{s1is1} and \ref{s1is2overk}}
\label{appC}

This appendix gives the proofs of Lemmas \ref{s1is2}, \ref{s1is1} and
\ref{s1is2overk}. The work here is basically a calculus marathon: while the
existence of solutions in the stated cases is relatively straightforward, the
nonexistence and uniqueness results are much more subtle.

We are looking for common zeros of the functions
\begin{align*}
h_1(x_1,x_2)&=\frac{2}{15}
(5x_2-5x_1+10s_1x_1^{2}-7x_1^{2}x_2-5x_1^{3}+2s_1x_1^{3}x_2)\\
h_2(x_1,x_2)&=\frac{2}{15}
(5x_1-5x_2+10s_2x_2^{2}-7x_2^{2}x_1-5x_2^{3}+2s_2x_2^{3}x_1).
\end{align*}
with $0<x_1<1$ and $0<|x_2|<1$ (where $x_2$ is negative if $g_2$ is negative
definite and positive if $g_2$ is positive definite). Since the equations
$h_1,h_2=0$ are both of the form $y(5-7x^2+2sx^3) - 5x + 10 s x^2 - 5 x^3=0$
we need to analyse the graphs of the functions $y= f_s(x): = \frac{5x (x^2 -
2s x + 1)}{2s x^3 - 7x^2 + 5}$ for $-1 < x < 1$.  Since also $|s_a|=2|{\mathbf
g}_a-1|/q_a$, where ${\mathbf g}_a$ is the genus of the corresponding curve
and $q_a\in \Z^+$, if $x_a$ is positive and $s_a>2/3$ then $s_a\in\{1,2\}$.
Thus for $s>2/3$ we can restrict out attention to the case where $-1 < x < 0$
or $s\in\{1,2\}$.  We then have the following lemma.

\begin{lemma}
\label{shapeofcurves}
Let ${\cC}={\cC}(s)$ denote the part of the graph of $y = f_s(x) = \frac{5x
(x^2 - 2s x + 1)}{2s x^3 - 7x^2 + 5}$ which lies within the square $[-1,1]
\times [-1,1]$. Then the following hold.

\begin{bulletlist}
  
\item When $0 \leq s \leq 2/3$,  ${\cC}$ looks like
\smallbreak
\centerline{\includegraphics[height=1in,width=1in]{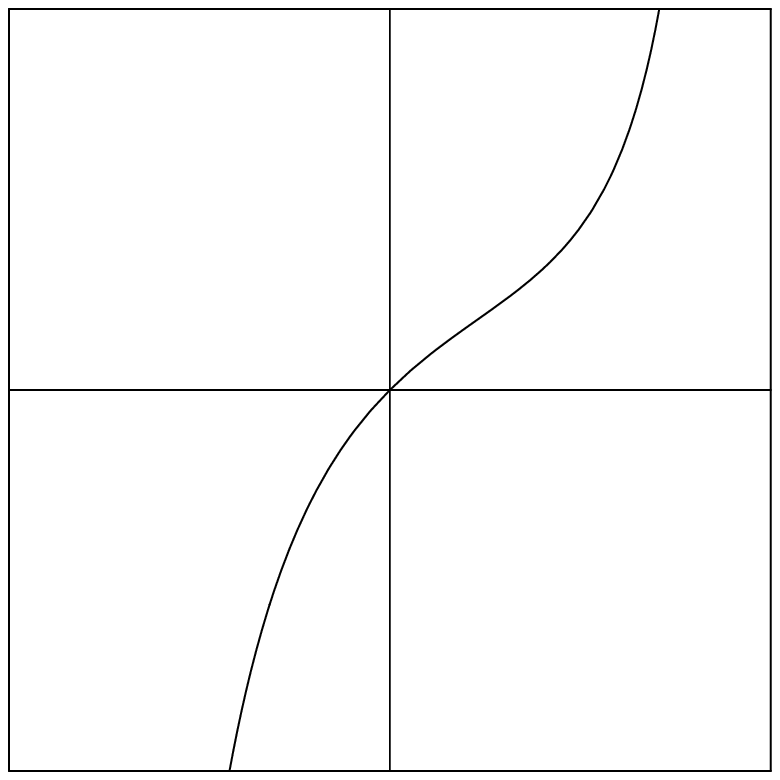}}
\noindent
where the graph is convex for $x<0$, increasing everywhere, intersects the
line $y=-1$ for some $-1<x<0$, and intersects $y=1$ for some $0<x<1$.

\item When $s=1$, ${\cC}$ looks like
\smallbreak
\centerline{\includegraphics[height=1in,width=1in]{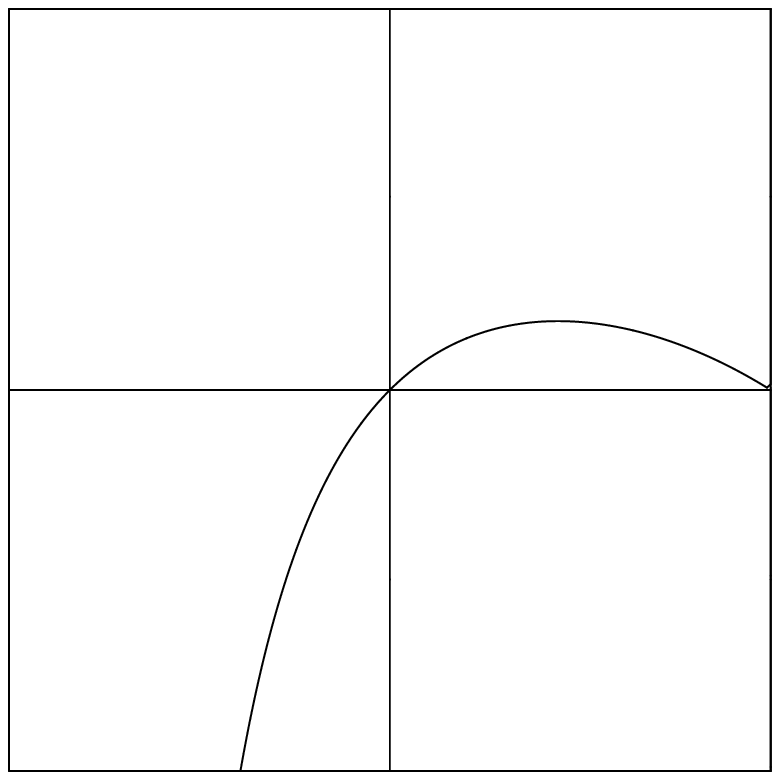}}
\noindent
where the graph is convex everywhere, increasing for $x<0$, intersects the
$x$-axis at $x=0$ and $x=1$, and intersects $y=-1$ for some $-1<x<0$.

\item  When $s=2$, ${\cC}$ looks like
\smallbreak
\centerline{\includegraphics[height=1in,width=1in]{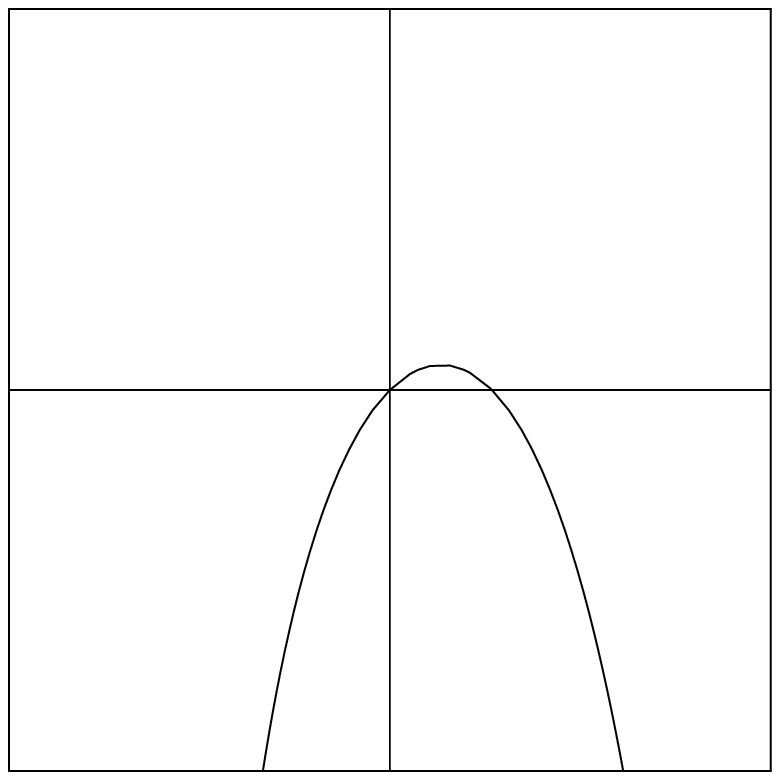}}
\noindent
where the graph is convex everywhere, increasing for $x<0$, intersects the
$x$-axis at $x=0$ and $x=2-\sqrt{3}$, and intersects $y=-1$ at $x=-1/3$ and
$x=(5 - \sqrt{10})/3$.

\item When $s \in (2/3, +\infty)$, ${\cC}$ restricted to $-1 < x < 0 $
looks like
\smallbreak
\centerline{
\includegraphics[height=1in,width=.5in]{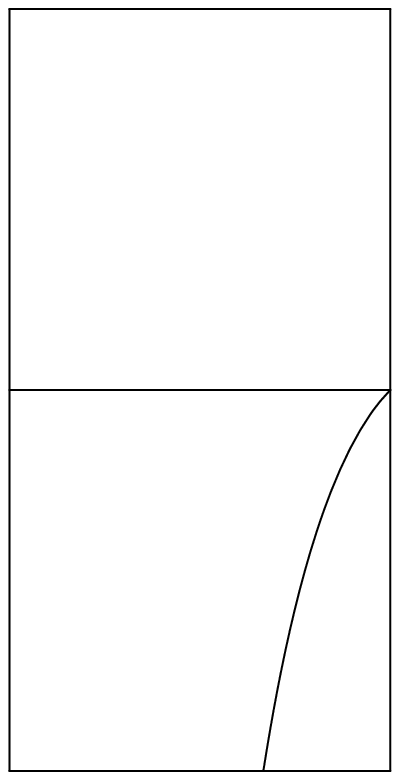}
\hbox to.5in{}}
\noindent
where the graph is convex and increasing, and intersects $y=-1$ for some $-1 <
x <0$.
\end{bulletlist}
\noindent
Since $-f_{-s}(-x)=f_s(x)$, for $s<0$, ${\cC}(s)$ is obtained by rotating
${\cC}(-s)$ by $\pi$.
\end{lemma}
\begin{proof}
The cases $s=0$, $s=1$ and $s=2$ are elementary and will be omitted.

We first consider the graphs for $-1<x<0$. The numerator of $\frac{5x (x^2 -
2s x + 1)}{2s x^3 - 7x^2 + 5}$ is strictly negative for $-1 < x < 0$, whereas
the denominator is negative at $x=-1$, positive at $x=0$ and strictly
increasing for $-1<x<0$. We conclude that $f_s$ has precisely one asymptote
$-1 < a < 0$ and $\lim_{x \rightarrow a^{\pm}} f_s(x) = \mp \infty$. Also
\begin{equation*}
f_s'(x) = \frac{5(5-20 s x + 22 x^{2} - 4s x^{3} - 7 x^{4} + 4 
s^{2}x^{4})}{(5 - 7 x^{2} + 2s x^{3})^{2}}
\end{equation*}
is positive for $-1 < x <0$ and since $f_s(-1) = 5 > 1$ the graph of $f_s$ is
outside the square $[-1,1] \times [-1,1]$ for $-1 < x < a$.  Thus we may
restrict our attention to $a < x < 0$ (and $f_s(x)$ is negative in this
range). Now
\begin{equation*}
f_s''(x) = \tfrac{-20(25 s - 90 x + 135 s x^{2} - 42 x^{3} - 70 s^{2} 
x^{3} + 51 s x^{4} - 6 s^{2} x^{5} - 7s x^{6} + 4s^{3} x^{6})}
{(5 - 7 x^{2} + 2s x^{3})^{3}}
\end{equation*}
is negative for $a < x < 0$ so $f_s$ is convex for $a < x < 0$. Since $\lim_{x
\rightarrow a^{+}} f_s(x) = -\infty$ and $f_s(0)=0$, ${\cC}$ must intersect
the line $y=-1$ for some $-1<a < x < 0$.

It remains to consider $0<x<1$ and $0 < s \leq 2/3$.

The denominator of $\frac{5x (x^2 - 2s x + 1)}{2s x^3 - 7x^2 + 5}$ is a third
degree polynomial which is negative at $x= -1$, positive at $x=0$, negative at
$x=1$ and positive for $x \rightarrow + \infty$, while the numerator is
positive for $0 < x < 1$. We conclude that $f_s$ has precisely one asymptote
$0 < b < 1$, $\lim_{x \rightarrow b^{\pm}} f_s(x) = \mp \infty$, $f_s(x) > 0$
for $0 < x <b$, and $f_s(x) < 0$ for $b < x < 1$.  For $x \in [0,1] \setminus
\{ b \}$ the denominator of
\begin{equation*}
f_s'(x) = \frac{5(5-20 s x + 22 x^{2} - 4s x^{3} - 7 x^{4} + 4 
s^{2}x^{4})}{(5 - 7 x^{2} + 2s x^{3})^{2}}
\end{equation*}
is positive. For a fixed $0 < x <1$, the numerator may be viewed as a function
of $s$ and its derivative, $5x(8 s x^{3} - 4 x^{2} - 20)$, with respect to
$s$ is clearly negative. Since the value of the numerator of $f_s'(x)$ at
$s=2/3$ equals
\begin{multline*}
\tfrac{5}{9}(45 - 120 x + 198 x^{2} - 24 x^{3} - 47 x^{4}) \\
=\tfrac{5}{9}(47(x^{2}-x^{4}) + 24(x^{2}-x^{3}) + (127 x^{2} - 120 x + 45)),
\end{multline*}
which is positive, we conclude that if $0 < s \leq 2/3$, $f_s'(x)$ is
positive for $x \in [0,1] \setminus \{b\}$.  Thus $f_s$ is strictly
increasing. Since $f_s(1) = -5$ the graph of $f_s$ is outside the square
$[-1,1] \times [-1,1]$ for $b < x < 1$.  Moreover, since $f_s(0)=0$ and
$\lim_{x \rightarrow b^{-}} f_s(x) = + \infty$, ${\cC}$ intersects the line
$y=1$ for some $0 < x < b<1$.
\end{proof}

It is clear from the shape of the graphs $\cC(s)$ (corresponding to $h_1=0$)
and their reflections in the line $y=x$ (corresponding to $h_2=0$) that the
zero-sets of $h_1$ and $h_2$ intersect in the fourth quadrant $0<x_1<1$,
$-1<x_2<0$ iff $s_1=2$ and $s_2=-2$, and in this case they meet at a unique
point $x_1=1/2$, $x_2=-1/2$. Hence we may assume from now on that $0<x_2<1$
and $s_2\leq 2$.

Let us now recall what we know about the functions $h_1$ and $h_2$:
\begin{itemize}
\item the curves $h_1=0$ and $h_2=0$ both pass through $(0,0)$;
\item along $h_1=0$ and $h_2=0$ we have $dx_2/dx_1=1$ at $(0,0)$;
\item along $h_1=0$ we have $d^2x_2/dx_1^2=-4s_1$ at $(0,0)$;
\item along $h_2=0$ we have $d^2x_2/dx_1^2=4s_2$ at $(0,0)$.
\end{itemize}
Therefore if $s_2>-s_1$ the zero-set of $h_2$ is above the zero-set of $h_1$
for $x_1$ small and positive, while if $s_2<-s_1$ it is below the zero-set of
$h_1$ for $x_1$ small and positive.

By Lemma~\ref{shapeofcurves}, the zero-sets of $h_2$ in $(0,1)\times(0,1)$ look
like
\smallbreak
\centerline{\includegraphics[height=.5in,width=.5in]{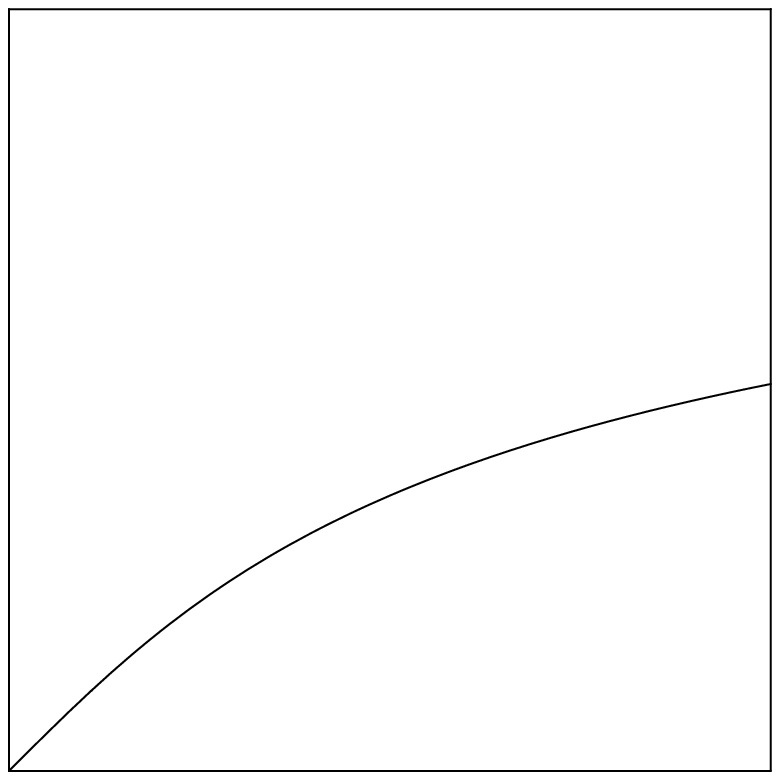}\qquad
\includegraphics[height=.5in,width=.5in]{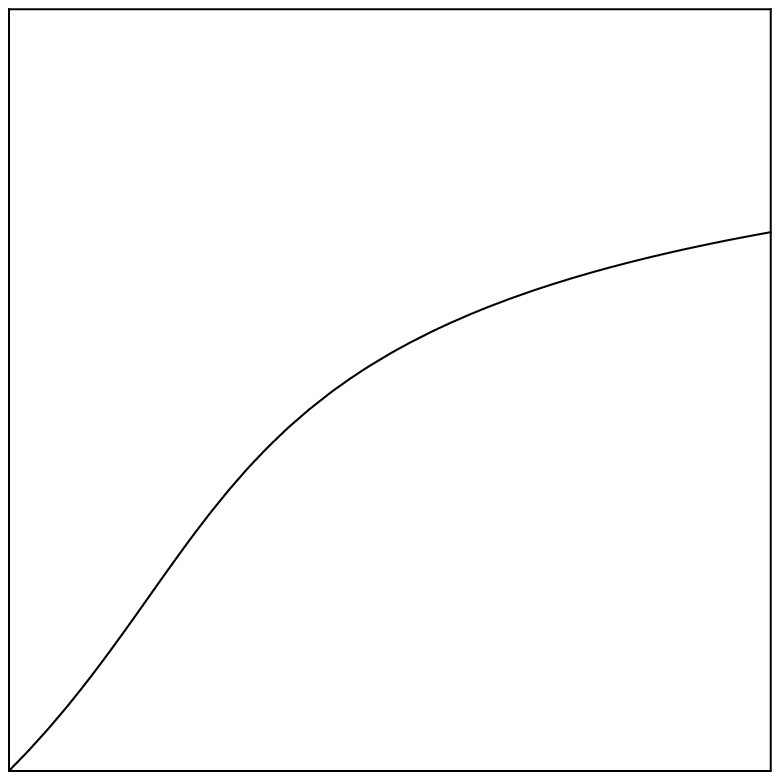}\qquad
\includegraphics[height=.5in,width=.5in]{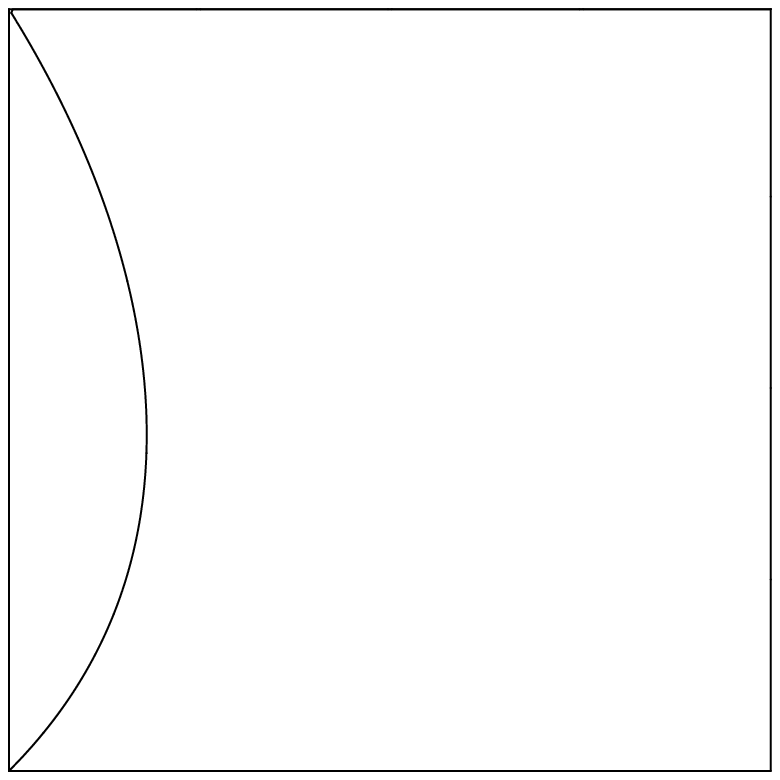}\qquad
\includegraphics[height=.5in,width=.5in]{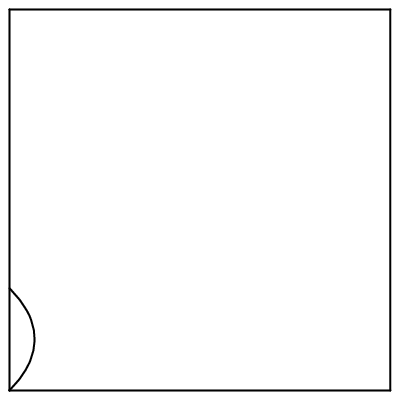}}
\noindent
for $s_2\leq 0$, $0<s_2\leq 2/3$, $s_2=1$ and $s_2=2$ respectively. For
$s_2\leq 2/3$ the zero-set of $h_2$ is an increasing graph which meets $x_1=1$
at a point with $0<x_2<1$.

It now follows easily that the zero-sets of $h_1$ and $h_2$ meet in at least
one point $(x_1,x_2)\in(0,1)\times (0,1)$ in the following cases:
\begin{itemize}
\item $s_1\in\{1,2\}$, $s_2<-s_1$;
\item $0<s_1\leq 2/3$, $-s_1<s_2\leq 2/3$.
\end{itemize}

For the nonexistence and uniqueness results, assume that we \emph{do} have a
solution $(x_{1},x_{2})\in (0,1) \times (0,1)$ for $h_{1}=h_{2}=0$. Then from
$h_{1}=0$, we have that
\begin{equation*}
x_{2}=\frac{5x_{1}(x_{1}^{2}-2 s_{1}x_{1} +1)}{2 s_{1} x_{1}^{3} - 
7x_{1}^{2} + 5}.
\end{equation*}
(It is easy to check that if $2 s_{1} x_{1}^{3} - 7x_{1}^{2} + 5 = 0$ then we
cannot have $h_{1}=0$ for $x_1\in (0,1)$.)  If we substitute into $h_{2}=0$ we
get
\begin{equation*}
\frac{10 x_1^2(x_1^2-5)\cM(x_1,1-x_1, s_1,s_2)}
{(2 s_1 x_1^3 - 7x_1^2 + 5)^3} = 0
\end{equation*}
with
\begin{align*}
\cM(x,y,s_1,s_2)&=-25 s_1 y^6 + 30 (2-5 s_1)xy^5 +20(15-23s_1)x^2y^4\\
&\quad+ 8(72-105s_1+25s_1^2)x^3y^3 +4(1-s_1)(132-125s_1+25s_1^2)x^4y^2\\
&\quad+8(1-s_1)^2(36-25s_1)x^5y + 96(1-s_1)^3x^6\\
&\quad- 25s_2 y(2x+y)\bigl(y^2+2(1-s_1)x(x+y)\bigr)^2.
\end{align*}
Thus if $(x_{1}, x_{2}) \in (0,1) \times (0,1)$ is any solution of
$h_{1}=h_{2}=0$, then $x_{1}$ must be a root of ${\cM}(x_1,1-x_1,s_{1},s_{2})$
and this determines $x_2$ uniquely. In the following $\cM_x$ denotes the
difference between the $x$ and $y$ derivatives of $\cM$, so
$\cM_x(x,1-x,s_1,s_2)$ is the $x$ derivative of $\cM(x,1-x,s_1,s_2)$;
$\cM_{xx}$ is defined similarly.

\subsubsection*{Proof of Lemma \textup{\ref{s1is2}}}
In this case $s_1=2$. We have seen that the zero-sets of $h_1$ and $h_2$ meet
in $(0,1)\times(-1,0)$ iff $s_2=-2$ and then the intersection point is unique,
being $(1/2,-1/2)$. We now analyse the case $x_2>0$.  Any intersection point
$(x_1,x_2)\in(0,1)\times(0,1)$ of the zero-sets of $h_1$ and $h_2$ must have
$0<x_{1}< 2- \sqrt{3}$ by Lemma \ref{shapeofcurves}, and $x=x_1$ must be a root
of $\cM(x,1-x,2,s_2)$ where
\begin{align*}
\cM(x,y,2,s_2)&= -96 x^6 - 112 x^5 y + 72 x^4 y^2 - 304 x^3 y^3 - 
    620 x^2 y^4 - 240 x y^5 - 50 y^6 \\
&\qquad    -25 s_2 y (2 x + y) (y^2 - 2 x (x + y))^2
\end{align*}
Clearly $\cM(x,1-x,2,s_2)$ is a decreasing function of $s_2$ when
$0<x<2-\sqrt{3}$. Since
\begin{equation*}
\cM(x,z+x,2,-2) 
= -4 x z (765 x^4 + 2040 x^3 z + 1846 x^2 z^2 + 680 x z^3 + 85 z^4),
\end{equation*}
$\cM(x,1-x,2,-2)<0$ for $0<x<2-\sqrt{3}<1/2$, hence so is
$\cM(x,1-x,2,s_2)$ for $s_2\geq -2$. Thus there are no solutions to
$h_{1}=h_{2}=0$ in $(0,1)\times(0,1)$ for $s_2\geq -2$.

Now suppose $s_{2}<-2$. We have seen that the zero-sets intersect in at least
one point $(x_{1},x_{2}) \in (0,1) \times (0,1)$. We now compute
\begin{align*}
\frac{\partial\cM_{xx}}{\partial s_2}(x,z+x,2,s_2)&=
-25 (9 x^4 + 216 x^3 z + 306 x^2 z^2 + 136 x z^3 + 17 z^4)\\
\cM_{xx}(x,z+2x,2,-2)&=
8 (4554 x^4 + 9340 x^3 z + 5757 x^2 z^2 + 1311 x z^3 + 85 z^4),
\end{align*}
so $\cM_{xx}(x,1-x,2,s_2)$ is a decreasing function of $s_2$ for
$0<x<2-\sqrt{3}<1/2$ whose value at $s_2=-2$ is positive for
$0<x<2-\sqrt{3}<1/3$. Hence $\cM(x,1-x,2,s_2)$ is a concave function of
$x\in(0,2-\sqrt{3})$. At $x=0$ it equals $-25(2 + s_2) >0$, while at
$x=2-\sqrt{3}$ it equals $48(240 - 139 \sqrt{3}) < 0$.  Hence it has exactly
one root $x=x_1\in(0,2-\sqrt{3})$ and the solution to $h_{1}=h_{2}=0$ is
unique.  \qed

\subsubsection*{Proof of Lemma \textup{\ref{s1is1}}}
In this case $s_{1}=1$. We have seen that the zero-sets of $h_1$ and
$h_2$ do not meet in $(0,1)\times(-1,0)$, so we restrict attention to
$x_2>0$. Since
\begin{equation*}
\cM(x,y,1,s_2) = -y^3 (64 x^3 + 160 x^2 y + 10 (9 + 5s_2) x y^2 +
25(1 + s_2)y^3)
\end{equation*}
there are no roots of $\cM(x,1-x,1,s_2)$ in $(0,1)$ for $s_2\geq-1$.  Suppose
now that $s_2<-1$. We have seen that the zero-sets intersect in at least one
point in $(0,1) \times (0,1)$.  The difference between the $x$ and $y$
derivatives of $-\cM(x,y,1,s_{2})/y^3$ is
\begin{equation*}
32 x^2 + 140 x y + 15 y^2 - 25 s_2 y (4 x + y)
\end{equation*}
which is clearly positive for $x\in(0,1)$, $y=1-x$ since $s_2$ is negative.
Hence $-\cM(x,1-x,1,s_{2})/(1-x)^{3}$ is an increasing function of $x\in
(0,1)$ so it has at most one root and the solution to $h_{1}=h_{2}=0$ is
unique.  \qed

\subsubsection*{Proof of Lemma \textup{\ref{s1is2overk}}}
In this case $s_{1} = 2/q_{1}$, $q_{1} = 3, 4, 5,\ldots$.  We have seen that
the zero-sets of $h_1$ and $h_2$ do not meet in $(0,1)\times(-1,0)$.  Thus we
may assume $0<x_2<1$ and $s_2\leq 2/3$: by the previous two lemmas (with
$s_1,s_2$ interchanged) there are no solutions with $s_2\in\{1,2\}$.

If there were a solution $(x_{1},x_{2})\in (0,1)\times(0,1)$ to
$h_{1}=h_{2}=0$ it would give a root $x=x_{1}$ of the function
${\cM}(x,1-x,s_{1},s_{2})$. We now observe that $\partial\cM/\partial s_{2}$
is negative for $y=1-x$, $x\in(0,1)$ (since $s_1<1$) and that
\begin{align*}
{\cM}(x,y,s,-s)&=4 x ((1 - s) x + y) \cM_0(x,y,s)\\
\notag
\cM_0(x,y,s)&=24 (1 - s)^2 x^4 + 48 (1 - s) x^3 y\\
\notag
&\qquad + 4 (21 - 25 s^2) x^2 y^2 + 20 (3 - 5 s^2) x y^3 + 5 (3 - 5 s^2) y^4
\end{align*}
so ${\cM}(x,1-x,s,-s)$ is positive on $(0,1)$ for $s\leq 2/3 <\sqrt{3/5}$.
Thus $\cM(x,1-x,s_1, s_{2})$ is positive for $0 < x < 1$ and $s_{2} \leq
-s_{1}$ and there are no solutions $(x_{1},x_{2})$ to $h_{1}=h_{2}=0$ with
$0<x_{1}< 1$ when $s_{2} \leq - s_{1}$.

We now let $s_2>-s_1$.  We have seen that the zero-sets intersect in at
least one point in $(0,1) \times (0,1)$. We want to show that they intersect
in at most one point. The proof, which is harder than previously, is motivated
by the observation that
\begin{equation*}
\cM(x,y,2/3,2/3)
=\tfrac{4}{27} (3 x^2 - x y - 5 y^2) (8 x^4 + 8 x^3 y + 112 x^2 y^2 +
120 x y^3 + 45 y^4)
\end{equation*}
and hence $\cM(x,1-x,2/3,2/3)$ is positive for $x_0<x<1$ where
$x_0=(9-\sqrt{61})/2$ is the smallest root of $x^2 - 9x + 5=0$.  Observe that
$x_0\approx 0.595$ is less that $3/5$ (since $1521=39^2$ is less than
$1525=5^2\cdot 61$). We are going to prove that $\cM(x,1-x,s_1,s_2)>0$ for
$3/5\leq x<1$ and that $\cM_x(x,1-x,s_1,s_2)>0$ for $0<x\leq 3/5$.  This will
prove that there is at most one root on $(0,1)$.

Since $\cM(x,1-x,s_1,s_2)$ is a decreasing function of $s_2$, to prove
positivity for $3/5\leq x<1$, it suffices to prove $\cM(x,1-x,s_1,2/3)>0$ for
$3/5\leq x<1$. This is true for $s_1=2/3$ and so the result follows from
\begin{claim} \label{c1}
$\frac{\partial\cM}{\partial s_1}(x,1-x,s_1,2/3)<0$ for $3/5\leq x<1$.
\end{claim}
The positivity of $\cM_x(x,1-x,s_1,s_2)$ on $0<x\leq 3/5$ for $-s_1<s_2\leq
2/3$ follows from the fact that it is an affine linear function of $s_2$ such
that:
\begin{claim} \label{c2}
$\cM_x(x,1-x,s_1,2/3)>0$ for $0<x\leq 3/5$;
\end{claim}
\begin{claim} \label{c3}
$\cM_x(x,1-x,s_1,-s_1)>0$ for $0<x\leq 3/5$.
\end{claim}
Subject to these three claims, we are done. \qed

\begin{proof}[Proof of Claim \textup1] We compute that
$-12\frac{\partial\cM}{\partial s_1}(z/2 + 3y/2,y,s,2/3)$ is given by
\begin{multline*}
18 (6665 - 11290 s + 6237 s^2) y^6 + (195155 - 361344 s +  194157 s^2) y^5 z\\
+ 2 (67267-131470 s +69255s^2) y^4z^2+2 (24931 - 50180 s + 26055 s^2) y^3z^3\\
+ 2 (5213 - 10610 s + 5445 s^2) y^2 z^4 + (1 - s) (1163 - 1197 s) y z^5 + 
54 (1 - s)^2 z^6.
\end{multline*}
It suffices to show that the coefficient of each monomial in $y,z$ is positive
for $0<s\leq 2/3$ (put $y=1-x, z=5x-3$). The first five
quadratics in $s$ have no real roots and are positive at $s=0$, and for the
last two the result is clear.
\end{proof}

\begin{proof}[Proof of Claim \textup2] We compute that $\frac{729}{2}
\cM_x(x,z/3 + 2x/3,s,2/3)$
is given by
\begin{multline*}
4 (88050 - 255955 s + 293700 s^2 - 99063 s^3) x^5
+80 (2460 - 4277 s + 5862 s^2 - 2025 s^3) x^4 z\\
+4 (19530 - 3229 s + 10575 s^2 - 4050 s^3) x^3 z^2
+4 (4632 + 2425 s - 975 s^2) x^2 z^3\\
+5 (420 + 347 s - 120 s^2) x z^4 + 10 (9 + 10 s) z^5.
\end{multline*}
It suffices to show that the coefficient of each monomial in $x,z$ is positive
for $0<s\leq 2/3$ (put $z=3-5x$). For the two quadratics and the last
coefficient, this is clear. The remaining three cubics are positive multiples
of
\begin{gather*}
44025 (2 - 3 s)^3 + 140270 (2 - 3 s)^2 s + 
    240345 (2 - 3 s) s^2 + 251028 s^3\\
1230 (2 - 3 s)^3 + 6793 (2 - 3 s)^2 s + 
    19272 (2 - 3 s) s^2 + 21789 s^3\\
9765 (2 - 3 s)^3 + 84656 (2 - 3 s)^2 s + 
    265431 (2 - 3 s) s^2 + 281844 s^3.
\end{gather*}
Hence they are all positive on $[0,2/3]$.
\end{proof}

\begin{proof}[Proof of Claim \textup3] We compute that $\frac{729}{4}
\cM_x(x,z/3 + 2x/3,s,-s)$
is given by
\begin{multline*}
48 (7125 - 23940 s + 23225 s^2 - 4974 s^3) x^5
+ 240 (1020 - 2259 s + 1054 s^2 + 525 s^3) x^4 z\\
+ 24 (4080 - 4509 s - 2750 s^2 + 4275 s^3) x^3 z^2\\
+  24 (897 - 450 s - 1225 s^2 + 750 s^3) x^2 z^3
+ 30 (25 - 6 s) (3 - 5 s^2) x z^4 + 30 (3 - 5 s^2) z^5.
\end{multline*}
It suffices to show that the coefficient of each monomial in $x,z$ is positive
for $0<s\leq 2/3$. This is clear for the last two
coefficients. The remaining four cubics are positive multiples of
\begin{gather*}
7125 (2 - 3 s)^3 + 16245 (2 - 3 s)^2 s - 
    2005 (2 - 3 s) s^2 + 363 s^3\\
510 (2 - 3 s)^3 + 2331 (2 - 3 s)^2 s + 
    2324 (2 - 3 s) s^2 + 1863 s^3\\
2040 (2 - 3 s)^3 + 13851 (2 - 3 s)^2 s + 
    22526 (2 - 3 s) s^2 + 15099 s^3\\
897 (2 - 3 s)^3 + 7173 (2 - 3 s)^2 s + 
    13919 (2 - 3 s) s^2 + 7419 s^3.
\end{gather*}
Only the first is not manifestly positive on $[0,2/3]$. However it is
positive at the endpoints and (dividing by $8$) the cubic $7125 - 23940 s +
23225 s^2 - 4974 s^3$ has a minimum at $s=(23225 - \sqrt{182167945}) /
14922\approx 0.652$ where it takes the value $5 (492445959775 - 36433589
\sqrt{182167945}) / 333999126\approx 10.5$, which is positive (since
$492445959775^2 = 242503023298720918050625>(36433589
\sqrt{182167945})^2 = 241810897419701928577345$).
\end{proof}

\begin{rem} The calculations in this final claim are remarkably tight. 
Numerical computations show that if we had broken the interval $(0,1)$ at a
point $\gtrsim 0.602$, instead of $3/5$, then this argument would fail, so we
are fortunate that $(9-\sqrt{61})/2\approx 0.595$ is less than this. 
We also remark that uniqueness of solutions to these equations can fail if we
allow $s_1,s_2\in(2/3,1)$, so the integrality conditions are crucial.

Depending on one's point of view, there are two possible responses to this
serendipity. The first is that it is just a coincidence that we obtain
unique WBF metrics in this (low-dimensional) situation. The second is that
there is a general uniqueness theorem for WBF metrics. We leave it to the
reader to decide.
\end{rem}

\end{document}